\documentclass[opre,nonblindrev,dvipsnames]{informs3}

\DoubleSpacedXI % Made default 4/4/2014 at request
\usepackage{endnotes}
\let\footnote=\endnote

\usepackage{natbib}
 \bibpunct[, ]{(}{)}{,}{a}{}{,}%
 \def\bibfont{\small}%

\usepackage{soul}

\TheoremsNumberedThrough     % Preferred (Theorem 1, Lemma 1, Theorem 2)
\ECRepeatTheorems

\EquationsNumberedThrough    % Default: (1), (2), ...
\usepackage{hyperref}
\usepackage{subfig}
\usepackage[ruled,vlined]{algorithm2e}
\usepackage{amsmath,amssymb,mathtools}

\newtheorem{lem}{Lemma}
\newtheorem{thm}{Theorem}

\newtheorem{cor}{Corollary}
\theoremstyle{definition}

\theoremstyle{remark}
\newtheorem{rem}{Remark}

\newtheorem{assum}{Assumption}

\newcommand{\E}{\mathbb{E}}
\newcommand{\Ea}[1]{\E\left[#1\right]}

\usepackage{changes}
\usepackage{color}
\usepackage{tikz}

\tikzset{
    ->,
    my_node/.style={
        align=center,
        minimum height=1.0cm,
        minimum width=2.2cm,
        rounded corners,
        draw,
    },
}

\definechangesauthor[color=blue]{TH}
\definechangesauthor[color=red]{AS}
\definechangesauthor[color=brown]{JG}

\begin{document}
%%%%%%%%%%%%%%%%

% \definecolor{r2}{HTML}{0000ff} % blue
\definecolor{r2}{HTML}{000000}
\definecolor{r3}{HTML}{000000}
\definecolor{r4}{HTML}{000000}
% \definecolor{r4}{HTML}{0000ff}

% Outcomment only when entries are known. Otherwise leave as is and
%   default values will be used.
%\setcounter{page}{1}
%\VOLUME{00}%
%\NO{0}%
%\MONTH{Xxxxx}% (month or a similar seasonal id)
%\YEAR{0000}% e.g., 2005
%\FIRSTPAGE{000}%
%\LASTPAGE{000}%
%\SHORTYEAR{00}% shortened year (two-digit)
%\ISSUE{0000} %
%\LONGFIRSTPAGE{0001} %
%\DOI{10.1287/xxxx.0000.0000}%

% Author's names for the running heads
% Sample depending on the number of authors;
% \RUNAUTHOR{Jones}
% \RUNAUTHOR{Jones and Wilson}
% \RUNAUTHOR{Jones, Miller, and Wilson}
% \RUNAUTHOR{Jones et al.} % for four or more authors
% Enter authors following the given pattern:
\RUNAUTHOR{Dias Garcia, Street, Homem-de-Mello, Muñoz}

\RUNTITLE{Application-Driven Learning Applied to Reserves and Demand Prediction}
% ARXIV/OO
% \RUNTITLE{Application-Driven Learning via Joint Prediction and Optimization of Demand and Reserves Requirement}

\TITLE{Application-Driven Learning: A Closed-Loop Prediction and Optimization Approach Applied to Dynamic Reserves and Demand Forecasting}

% Block of authors and their affiliations starts here:
% NOTE: Authors with same affiliation, if the order of authors allows,
%   should be entered in ONE field, separated by a comma.
%   \EMAIL field can be repeated if more than one author
\ARTICLEAUTHORS{%
\AUTHOR{Joaquim Dias Garcia}
\AFF{LAMPS, DEE, PUC-Rio \& PSR, Rio de Janeiro, Brazil, \EMAIL{joaquim@psr-inc.com}} %, \URL{}}
\AUTHOR{Alexandre Street}
\AFF{LAMPS, DEE, PUC-Rio, Rio de Janeiro, Brazil, \EMAIL{street@ele.puc-rio.br}}
\AUTHOR{Tito Homem-de-Mello}
\AFF{School of Business, Universidad Adolfo Ibáñez, Santiago, Chile, \EMAIL{tito.hmello@uai.cl}}
\AUTHOR{Francisco D. Muñoz}
\AFF{Generadoras de Chile, Santiago, Chile, \EMAIL{fmunoz@generadoras.cl}} %, \URL{}}
% Enter all authors
} % end of the block

\ABSTRACT{%
%200 words maximum

{\color{r3}
Decision-making is generally modeled as sequential forecast-decision steps with no feedback, following an open-loop approach. For instance, in the electricity sector, system operators use the forecast-decision approach followed by ad hoc rules to determine reserve requirements and biased net load forecasts to guard the system against renewable generation and demand uncertainty. Such procedures lack technical formalism to minimize operating and reliability costs.  We present a new closed-loop framework, named application-driven learning, in which the best forecasting model is defined according to a given application cost function. We consider applications in which the decision-making process is driven by two-stage optimization schemes fed by multivariate point forecasts. We present our estimation method as a bilevel optimization problem and prove convergence to the best estimator regarding the expected application cost. We propose two solution methods: an exact method based on the KKT conditions of the second-level problems and a scalable heuristic suitable for decomposition. Thus, we offer an alternative scientifically grounded approach to current ad hoc procedures implemented in industry practices. We test the proposed methodology with real data and large-scale systems with thousands of buses. Results show that the proposed methodology is scalable and consistently performs better than the standard open-loop approach.

}
}%

% Fill in data. If unknown, outcomment the field
\KEYWORDS{Application-driven learning, {\color{r2}closed-loop} prediction and optimization, bilevel optimization, {\color{r2}dynamic}
 reserves, forecast, power-systems {\color{r2}operation}}
%\HISTORY{This is a working paper. First ideas were drafted in Santiago in Nov 2019. Then discussion followed in Rio in Jan 2020.}

\maketitle
%%%%%%%%%%%%%%%%%%%%%%%%%%%%%%%%%%%%%%%%%%%%%%%%%%%%%%%%%%%%%%%%%%%%%%

% Samples of sectioning (and labeling) in OPRE
% NOTE: (1) \section and \subsection do NOT end with a period
%       (2) \subsubsection and lower need end punctuation
%       (3) capitalization is as shown (title style).
%
%\section{Introduction.}\label{intro} %%1.
%\subsection{Duality and the Classical EOQ Problem.}\label{class-EOQ} %% 1.1.
%\subsection{Outline.}\label{outline1} %% 1.2.
%\subsubsection{Cyclic Schedules for the General Deterministic SMDP.}
%  \label{cyclic-schedules} %% 1.2.1
%\section{Problem Description.}\label{problemdescription} %% 2.

% Text of your paper here

\section{Introduction}\label{sec-intro}

The most common approach to making decisions under uncertainty involves two steps.
In the first step, one develops a forecast for all uncertainties that affect the decision-making problem based on all information available. In the second step, an action {\color{r3} is planned and implemented based on the forecast. In the case of sequential decision-making problems, however, this decision step may also consider some corrective actions} after uncertainties are realized. The aforementioned two-step procedure constitutes an open-loop forecast-decision process in which the outcomes (costs) of the decisions {\color{r3}actually implemented} are not accounted for in the forecasting {\color{r3}step}.

In the electricity sector, it is common for system operators to use an open-loop forecast-decision approach. First, loads are forecast based on standard statistical techniques, such as least squares (LS), and reserve requirements are defined by \emph{ad hoc} rules based on quantiles, extreme values, and standard deviation estimates of forecast errors following specified reliability standards \citep{NRELguidelines}. Then, a decision is made to allocate generation resources following an energy and reserve scheduling program \citep{Chen2014, Tony2019}. In real-time, reserves are deployed to ensure power is balanced at every node, compensating for forecast errors. 

{\color{r3}There is relevant empirical evidence that system operators rely on \textit{ad hoc} or out-of-market actions---and not just on reserves---to deal with uncertainty and cost asymmetry in power system operations.
For instance, according to the 2019 Annual Report on Market Issues and Performance of the California ISO \citep{CAISO}, \textit{``...operators regularly take significant out-of-market actions to
address the net-load uncertainty over a longer multi-hour time horizon (e.g., 2 or 3 hours). These actions
include routine upward biasing of the hour-ahead and 15-minute load forecast, and exceptional
dispatches to commit and begin to ramp up additional gas-fired units in advance of the evening ramping
hours."} Additionally, reserve requirements are, in practice, empirically defined according to further \emph{ad hoc} off-line rules based on off-line analysis \citep{NRELguidelines,PJM2018}.} 

These \emph{ad hoc} forecast biasing procedures to minimize operating and reliability costs lack technical formalism and transparency, which prevents agents from internalizing them.
{\color{r3}However, such actions are based on the empirical evidence system operators acquire from real-world practice, i.e., from the acknowledgment of the fact that a positively biased demand forecast (above the expected value) may lead to cheaper corrective actions (reduce or curtail scheduled generation is much cheaper than increase it or shed load).} Consequently, {\color{r3}there is a potential and eminent benefit to be unlocked in real-world applications based on the traditional deterministic (forecast-decision) approach by closing the loop between the prediction and prescription steps. This is the main focus of this paper. 
\subsection{Gap in the literature}

{\color{r4}To accurately position our research within the existing gaps in the literature, we start by highlighting that when we refer to ``forecasts", we mean \textit{pointwise forecasts}. Additionally, these forecasts will be used in deterministic decision (planning) models.} One might correctly question why we do not consider more complete forecasts that yield a probability distribution instead of a single estimate. 
Indeed, even in the context of our application, it has been demonstrated that stochastic programming models yield better results than deterministic ones when making decisions under uncertainty because the former takes distributions into consideration. These models provide better results in terms of cost, reliability, and market efficiency compared to deterministic approaches \citep{wang2014flexible}. 

However, in practical applications, such as the case of {\color{r4}unit commitment and} scheduling energy and reserves, where multiple stakeholders are involved, the following issue arises: proper modeling of high-dimensional distributions in stochastic optimization models is challenging, and tractability requires the application of the well-known sample average approximations (SAA) approach. As a consequence of applying the SAA to solve practical stochastic optimization problems, solutions become sample dependent (see \cite{papavasiliou2014applying, papavasiliou2013multiarea}), thereby compromising market transparency and preventing stakeholders' acceptance \citep{wang2015real}. Therefore, most system operators worldwide still rely on deterministic short-term scheduling (economic dispatch or unit commitment) models with exogenous forecasts for loads and reserve requirements \citep{Chen2014, PJM2018}. 

{\color{r4}Other attempts to incorporate probabilistic information have been proposed in the literature. For instance, although time-varying reserves have been studied in the past, 
they have re-emerged as dynamic probabilistic reserves \citep{Tony2019}. 
In the context of a more sophisticated dynamic probabilistic reserve approach, probabilistic forecasts are frequently used to account for forecast errors. These probabilistic reserves can be sized following a variety of methods with different complexity based on forecast error standard deviations \citep{strbac2007impact, holttinen2012methodologies},
non-parametric estimation of the forecast error distribution \citep{bucksteeg2016impacts},
or even machine learning \citep{Tony2019}. These are all considered stochastic methods and are simple alternatives to capture and incorporate fairly complex dynamics that are challenging for bottom-up approaches.} {\color{r4} Nonetheless, all previously reported works can be seen as open-loop approaches, wherein more sophisticated reserve requirements are fed into deterministic models. So, despite providing relevant contributions to the subject, mainly by the dynamic adaptation of forecasted reserve levels based on probabilistic information, the application-cost feedback on demand forecasts and reserve requirements is still absent.}

{\color{r4}Given that} changing the decision-making framework (prescriptive model) from deterministic to stochastic might constitute a major and possibly impractical task, {\color{r4}it is no surprise to find evidence of system operators biasing their point-forecasting framework (predictive model) to obtain better results in real time} (\cite{CAISO}). \emph{Thus, any improvement in the latter has the potential to be directly absorbed by system operators and impact industry practices and outcomes.}  {\color{r4}As mentioned earlier, we do so by developing a closed-loop approach to determine an optimal pointwise forecast}.

{\color{r4}Taking a close look at the literature on the subject, }{\color{r3}the idea of ``closing the loop'' by adapting the forecast to the application at hand dates back at least to the 1980s. It has long been observed that }classical forecasting methods do not take the underlying application of the forecast into account. Consequently, hypotheses such as prediction error symmetry in the LS might not be the best fit for problems with asymmetric outcomes. By acknowledging the asymmetry in particular problems, researchers have attempted to capture it empirically; however, such an approach does not take the application into account directly. Some existing methods do capture asymmetry, such as Quantile Regression \citep{rockafellar2008risk}. Also, \cite{zellner1986bayesian} and \cite{zellner1986biased} acknowledge that biased estimators can perform even better than those that make accurate predictions of statistical properties of the stochastic variables. The authors exemplify that an overestimation is not as bad as an underestimation for the case of dam construction and attribute a second example about the asymmetry on real estate assessment to \cite{varian1975bayesian}. {\color{r3}In the machine learning community, which combines many ideas from optimization and probability, it was also proposed long ago to treat the prediction and prescription steps jointly \citep{bengio1997using}.}

{\color{r4} More recently, the explosion in the availability of data---including \textit{contextual information}–––and in the development of machine learning methods has led to relevant efforts by the operations research community to integrate the forecast and optimization steps; a very recent survey by \cite{sadana2023survey} describes the recent advances in this area.
There is a stream of literature that aims at said integration not by ``closing the loop" as discussed above but by changing the optimization step in order to fully benefit from the use of machine learning methods in the forecast step. Examples of such approaches, called ``predict-then-optimize" or ``estimate-then-optimize" in the literature, include \cite{bertsimas2019predictive}, \cite{ban2019dynamic}, \cite{DiaoSen2020} and \cite{kanBL:20}. On the other hand, new techniques have been studied to combine the forecast and optimization steps using the structure of the application at hand in the forecast step. Such an approach has received different names, such as ``smart predict-then-optimize" \citep{elmachtoub2017smart}, ``integrated estimation and optimization" \citep{Grigas-ICEO:2021}, ``end-to-end learning" \citep{donti2017task}, and ``decision-focused learning" \citep{mandi2022decision}. Those works, however, do not tackle the type of models we study in this paper, where we deal with polyhedral functions (as opposed to linear or smooth functions) and very large-scale problems.}

Within this context, two possible avenues of research are opened to achieve better results: i) a focus on improving the decision-making model (prescriptive framework), which assumes we can change it to incorporate better forecast information (e.g. forecasts that include contextual information); or ii) a focus on improving the forecasting model (predictive framework), which assumes we can not change the decision-making process (in our application, defined by system operator's dispatch model), but we can change the forecasts to incorporate, in a closed-loop manner, a given application cost function.
{\color{r3}The latter is the more realistic assumption for power system operations, as changing dispatch models requires a great technical and political effort, whereas biasing the load forecast and adjusting reserve requirements according to \emph{ad hoc} rules constitute the current practice in this sector worldwide.} {\color{r4}  As discussed below, we focus on the latter avenue. }

{\color{r2}

\subsection{Objective and contribution}\label{sec-con}

{\color{r3}The objective of this paper is to present a new closed-loop framework, named \textit{application-driven learning}, in which the best point forecasting model is defined according to a given application cost function that can be represented by a two-stage linear program with uncertainty on the right-hand side. Note that, because we are dealing with point forecasts, the two-stage model reduces to a linear program with uncertainty on the right-hand side. In the proposed method, the \emph{application schema} is characterized through an optimization model, which is then used as part of the estimation problem. Thus, our framework replaces the traditional statistical error minimization objective with a cost-minimization structure of a specific application. Although the method is general, in this paper, we study an application in electricity markets, as this was the main motivation for these developments and is one of the most relevant examples of industry applications in which the decision-making process is driven by a large-scale optimization scheme fed by multivariate point forecasts.

To achieve the objectives described before, we derive the following general technical contributions: 

\begin{enumerate}
    
   \item We develop a new and flexible application-driven learning framework to yield a point forecast that performs best in a given application based on a bilevel optimization model (Section \ref{sec-theory}). To the best of the authors' knowledge, for the first time in the literature, the relevant case of integrated estimation and optimization for applications based on linear programming models affected by right-hand-side uncertainty is addressed with specialized algorithms. Two solution approaches are presented in Section \ref{sec-solution}. The first approach is an exact method based on the KKT conditions of the second-level problem. The second is a scalable heuristic approach suitable for decomposition methods and parallel computing. Although not limited to linear bilevel programs, we show how to design efficient methods tailored for off-the-shelf linear optimization solvers. Additionally, a salient feature of our scalable heuristic method is that it ensures optimal second-level solutions, which emulates the reality where the plan is optimized in the application context.
    Hence, the proposed framework is general and suitable for a wide range of applications relying on the standard structure of the forecast decision process. 

    \item We provide asymptotic convergence proofs for both the objective function value and estimated parameters of the proposed application-driven learning method in Section \ref{sec-conv}. The convergence proof highlights that the method is asymptotically the best that can be done for a specific application (described by planning and implementation processes) given a forecast functional form (see Corollary \ref{cor}). 
    Under the hypothesis of our methodology and based on the aforementioned convergence results, we show that the solution of our method converges to the best estimator in terms of the expected cost of the selected application.

    \item In addition to the above methodological contributions, our paper also provides relevant contributions to the applied field of Operations Research (OR) in energy. In particular,
    we propose a new methodology to dynamically forecast the load and define the reserve requirements for the problem of scheduling energy and reserves in Section \ref{sec-power}. The method can be used to either jointly optimize the load forecast and reserve requirements or to optimize only the reserve requirements given an exogenous forecast for the load. In both cases, the optimal solution defines the optimal policy to dynamically allocate reserve requirements so that the expected dispatch cost is minimized in the long run. Based on this contribution, agents are provided with a scientifically grounded and comprehensively described methodology that can be used to reduce the number of \emph{ad hoc} procedures currently implemented in practice.

    \item {\color{r4}We present an additional methodology for load forecasting that mathematically describes the load biasing practice seen in the industry in Section \ref{sec-naive}. The method obtains the best multiplier for a previously estimated forecast model. Therefore, it is a naive \textit{ad hoc} linear biasing benchmark and, as expected, will lead to suboptimal solutions compared to the application-driven learning framework.}
\end{enumerate}

Finally, to empirically corroborate the relevance of our contributions and to demonstrate the applicability, performance, and scalability of our methodology, we benchmark the proposed method with the traditional sequential least squares forecast and energy and reserve scheduling approach (Section \ref{sec-simulations}). To do that, we analyze the proposed methodology in several case studies using multiple test systems ranging from single bus systems to realistic, very large systems with up to 13,659 buses. In our tests, we consider both synthetic and real load data to demonstrate the effectiveness of our proposed approach. 
Results show that the two proposed application-driven learning approaches (demand and reserves and only reserves) consistently yield better performance on out-of-sample tests than the benchmark, where forecasts and decisions are sequentially carried out. 
For very large systems, where the exact method fails to find solutions within reasonable computational times, the heuristic method exhibits high-quality performance compared to the benchmark for all test systems. 
{\color{r4}As benchmarks, we consider both LS-based forecasts and \textit{ad hoc} biased forecasts as mentioned in \cite{CAISO}.}
}

\section{Application-Driven Learning and Forecasting}\label{sec-theory}

In this section, we contrast the standard sequential framework, referred to as \textit{open-loop}, and the joint prediction and optimization model, referred to as \textit{closed-loop}. The presentation is in general form to facilitate the description of the solution algorithm, to set notation for the convergence results, and to highlight that the method has applications beyond load forecasting and reserve sizing in power systems. We will specialize the bilevel optimization problem for closed-loop load forecasting and reserve sizing in Section \ref{sec-power}.

We consider a dataset of historical data $\{y_t, x_t\}_{t\in \mathbb{T}}$, where $\mathbb{T} = \{1,\dots, T\}$. Here $y_t \in \mathbb{R}^n$ are observations of a variable of interest that we want to forecast, while $x_t \in \mathbb{R}^m$ are observations of external variables (covariates or features) that can be used to explain the former. Furthermore, the latter might include lags of $y_t$ as in autoregressive time series models.

The classic forecast-decision approach works as follows. The practitioner \textit{trains} a parametric \textit{forecast} model seeking for the best vector of parameters, $\theta$, such that a loss function, $l(\cdot,\cdot)$, between the conditional forecast for sample $t$, $\hat{y}_{t}(\theta, x_t)$, and the actual data, $y_{t}$, is minimized, i.e., solving $\min_\theta \frac{1}{T} \sum_t l( \hat{y}_{t}(\theta, x_t), y_{t})$. This is frequently done by solving LS optimization problems and finding $\theta^{LS} \in \argmin_\theta \frac{1}{T} \sum_t \| \hat{y}_{t}(\theta, x_t) - y_{t}\|^2$.
In the \textit{planning} step, a decision is made by an optimized policy based on the previously obtained forecast, i.e., with $\hat{y}_{t}^{LS} = \hat{y}_{t}(\theta^{LS}, x_t)$. This results in a vector $z^*(\hat{y}_{t}^{LS})$, which in our application comprises the schedule of energy and reserves through generating units. Finally, the actual data $y_t$ is observed, and the decision-maker must adapt to it (for instance, the system operator may respond with a balancing re-dispatch to adjust the generation within the scheduled reserves), and a \textit{cost}, $G_a( z^*(\hat{y}_{t}^{LS}), y_{t})$, is measured (for instance, based on the final adjusted generation cost). There is no feedback from the final cost in the forecasting and decision policy, hence, the name \textit{open-loop}.

\subsection{The proposed closed-loop application-driven framework}\label{sec-closed}

The core of the proposed predictive framework is to explore a feedback structure between the estimated predictive model and the application cost assessment. The general idea is depicted in Figure \ref{fig-closedloop}, which also stresses the difference from the open-loop model.

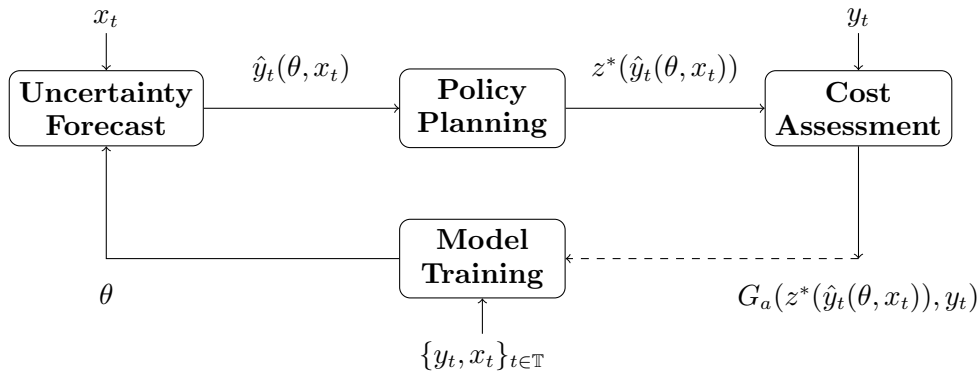
\begin{figure}[h]
\centering
\begin{tikzpicture}
    \linespread{0.5}
\node[my_node] at (-5,0) (forecast) {\textbf{Uncertainty} \\ \textbf{Forecast}};
\node[my_node] at (-0,0) (planning) {\textbf{Policy}\\ \textbf{Planning}};
\node[my_node] at (+5,0) (assessment) {\textbf{Cost}\\ \textbf{Assessment}};
\node[my_node] at (0,-2) (train) {\textbf{Model}\\  \textbf{Training}};
\coordinate (assess_train) at (-5,-2); 
\coordinate (train_forecast) at (+5,-2); 
\coordinate (dataset) at (0,-3); 
\coordinate (open_loop_end) at (+5,-2);
\coordinate (x) at (-5,+1); 
\coordinate (y) at (+5,+1); 
\draw[->] (x) -- node [above=0.5em,]{ $x_t$} (forecast);
\draw[->] (y) -- node [above=0.5em,]{ $y_t$} (assessment);
\draw[->] (forecast) -- node [above=0.5em]{ $\hat{y}_t({\theta}, x_t)$} (planning);
\draw[->] (planning) -- node [above=0.5em]{$z^*(\hat{y}_t(\theta, x_t))$} (assessment);
\draw[->] (assessment) -- node [below=1.5em, left=3em]{} (open_loop_end);
\draw[dashed, ->] (open_loop_end) |- node [below=0.5em]{$G_a(z^*(\hat{y}_t(\theta, x_t)), y_t)$} (train);
\draw[->] (train) -| node [below=0.5em]{${\theta}$} (forecast);
\draw[->] (dataset) -- node [below=0.5em,]{$\{y_t, x_t\}_{t\in\mathbb{T}}$} (train);
\end{tikzpicture}
\caption{Learning models: considering the dashed line we have the \textit{closed-loop} model, otherwise it represents the \textit{open-loop} model.} \label{fig-closedloop}
\end{figure}

The estimation method can be mathematically described through the following bilevel optimization problem: 
\begin{align}
    \theta_T {\color{r2}\in} \argmin \limits_{\theta \in \Theta, \hat{y}_t,z_t^{*}} \quad&  \frac{1}{T} \sum_{t \in \mathbb{T}} G_a(z_t^{*}, {y}_t) \label{mod-main-obj} \\
    s.t. \quad& \hat{y}_t = \Psi (\theta, x_t) \ \forall t \in \mathbb{T} \label{mod-main-forecast}\\
        & z_t^{*} \in \argmin_{z \in Z }  G_p(z, \hat{y}_t)\ \forall t \in \mathbb{T}, \label{mod-main-policy}
\end{align}\label{model}
In model \eqref{mod-main-obj}--\eqref{mod-main-policy}, $\Psi (\theta, x_t)$ represents a \textit{forecasting} model that depends on both the vector of parameters, $\theta$, and the features vector, $x_{t}$, possibly including lags of $y_t$. The vector $\hat{y}_t$ is the forecast generated for sample (or period) $t$ conditioned to the vector of features, ${x}_{t}$, as defined in \eqref{mod-main-forecast}. For each $t$, the forecast $\hat{y}_t$ is used as input in a lower-level problem and a planning policy, $z_t^{*}$, is obtained as a function of $\hat{y}_t$, i.e., $z_t^{*}(\hat{y}_t)$. This is done by optimizing the decision-maker \textit{planning} cost function, $ G_p(z, \hat{y}_t)$ in \eqref{mod-main-policy}. Then, the optimized policy $z_t^{*}$  is evaluated in the first level against the actual realization, $y_t$, for each $t$. The evaluation is made under the decision-maker's \textit{assessment} (or \textit{implementation}) cost function, $G_a(z_t^{*}(\hat{y}_t), {y}_t)$. Hence, the application is embedded into the estimation process in both the \emph{ex-ante} planning policy and \emph{ex-post} assessment objective \eqref{mod-main-obj}--\eqref{mod-main-policy}.

It is worth noticing that the proposed formulation can be interpreted as an optimization over $\theta$ in a back-test, in which for a given $\theta$, the assessment of the forecast performance is completely determined by $G_a(z_t^{*}(\hat{y}_t), {y}_t)$. Within this context, the upper level identifies the parameters with the best back-test performance.
Furthermore, note that for a fixed $\theta$, there is no coupling between two samples. Therefore, the model can be decomposed per $t$. This will be used in one of the proposed solution methods. 
Finally, there is a slight abuse of notation in \eqref{mod-main-obj} because the argmin only retrieves $\theta_T$, {\color{r2}a} solution for $\theta$ with $T$ samples, thus disregarding the rest of the upper level decision vectors, $\hat{y}_t$ and $z_t^{*}$. {\color{r2}Note that $\theta_T$ is one among the multiple options in the set of all possible solutions $S_T$}.

One key difference from previous works \citep{donti2017task, elmachtoub2017smart, munoz2022bilevel} is that $G_a$ and $G_p$ can be different functions. 
{\color{r4} As we shall see later, in the case of our application, $G_p$ corresponds to the system operator's description of the system considered in the dispatch optimization, i.e. the \textit{planning phase}, whereas $G_a$ corresponds to an \textit{assessment phase}, which evaluates the plan given by the solution to $G_p$.}
This is extremely useful in the context of power systems operations where planning models might differ from real-time ones. {\color{r4} For instance, while in some market designs, zonal models are considered in the planning phase to simplify transmission constraints, a nodal description can still be used to evaluate the implemented dispatch.}

Next, we specialize the above general model to a linear one, so, for $i\in \{a, p\}$ we have:
\begin{align}
    & G_i(z, y) = c_i^{\top}z + Q_i(z, y) \label{mod-main-mod} \\
    & Q_i(z, y) = \min \limits_{u} \{q_i^{\top} u \ | \ W_iu \geq b_i - H_i z + F_i y \} \label{mod-main-rec}
\end{align} \label{recourse}
\noindent Note that the functions in \eqref{mod-main-mod} and \eqref{mod-main-rec} resemble the formulation of two-stage stochastic programs, in the sense that given a decision $z$ and an observation $y$, one determines the best corrective action $u$. In that context, $c_i,\ q_i,\ W_i,\ b_i,\ H_i$ and $F_i$ ($i \in \{a,p\}$) are parameters defined according to the problem of interest. Note also that the uncertainty $y$ appears only on the right-hand side of the problems defining $Q_a$ and $Q_p$; this will be important for our convergence analysis and solutions methods.

Although model \eqref{mod-main-obj}--\eqref{mod-main-policy} is fairly general, we specialized to the case of linear programs {\color{r4} with} right-hand-side uncertainty \eqref{mod-main-mod}--\eqref{mod-main-rec}, because we will assume polyhedral structure for the set $Z$. This can be contrasted with previous works {\color{r4}  such as \citep{donti2017task}, which considers strongly quadratic programs, and \citep{elmachtoub2017smart}, which considers linear programs with uncertainty in the objective function. In our case, the objective functions (in both upper and lower levels) correspond to the objective value of a linear program as a function of its right-hand side and as such are polyhedral functions. As mentioned earlier, the polyhedral structure}  will be important for developing our asymptotic convergence results and our solution methods.

{
\color{r2}
\section{A Motivating Example}

In this interlude, we present a small stylized example to showcase how asymmetries can affect open-loop and closed-loop models. Although the model is very simple, we use the optimization notation from the previous section to illustrate the ideas discussed there.

Consider the scheduling process for the next hour of a power system containing a single power plant with a capacity of $4$ MW (hence, capable of generating $4$ MWh in an hour) and a generation cost of \$10/MWh. Consider a penalty of \$100/MWh for scheduling the plant below the realized demand and \$0/MWh for scheduling above the realized demand. 
{\color{r3}The goal of the system operator is to determine the optimal amount of energy $z^*$ to be generated, given a pointwise forecast of the demand $\hat{y}$. That is, the operator solves the planning problem  given below:
\begin{align}
\label{mod:planning}
    \min_{z\in [0,4]} \ & \big\{G_p(z,\hat{y}):= 10 z+ 100 (\hat{y}-z)_+ + 0 (z-\hat{y})_+ \big\}.
\end{align}
After the decision is made, the actual incurred cost is given by the assessment function $G_a$ {\color{r4} (note that, in this example, the coefficients of $G_p$ and $G_a$ are the same for the sake of simplicity, but as discussed earlier they may be different in general):}
\begin{align}
\label{mod:assess}
     G_a(z^*,y):= 10z^* + 100 (y-z^*)_+ + 0 (z^*-y)_+ 
 \end{align}
where $z^*$ is the fixed decision from the planning problem, and $y$ is the actual demand that is realized during the implementation. }

Suppose that we directly forecast the value of demand, i.e., the forecasting model is $\hat{y} =\Psi(\theta) = \theta$, and that 50\% of data is equal to 0 and 50\% is equal to 2 MWh.  Consider first the situation in which the operator uses a standard LS estimator for the demand. Then, the demand forecast is $\hat{y} = 1$ MWh, {\color{r4}  and a simple calculation shows this corresponds to taking $z^*=1$, thus leading to an average cost of \$60. } 

{\color{r4} However, based on our methodology, we can obtain the \textit{best} forecast for this application.  As seen earlier,  in this example we have $\hat{y}= \theta$, hence the minimization problem in \eqref{mod-main-obj} can be expressed as 
\begin{equation}
\label{eq:toy}
  \theta_T \in \argmin \limits_{\theta} \, \frac{1}{2} \big(G_a(z^*, 0) + G_a(z^*, 2)\big),
\end{equation}
where $G_a$ is given in \eqref{mod:assess} and $z^*=z^*(\theta)$ is the optimal solution to \eqref{mod:planning} when $\hat{y}= \theta$. So, we see that the model becomes a \textit{bilevel} one.  In this ``toy" example, the lower level problem has an analytical solution $z^*(\theta)=\theta$, so simple calculations show that the minimizer of \eqref{eq:toy} is $\theta^* = 2$. This yields the \textit{optimally biased} forecast $\hat{y} = 2$ (MWh) and the corresponding planning solution $z^*=2$, which in turn leads to an average cost of just \$20.}

One of the main reasons for the phenomenon described above is that the operator is constrained to a specific decision rule (planning method) that has to consider a given forecast in a predetermined way. While the above example is indeed very stylized to allow for a simple exposition, it carries the fundamental idea that {\color{r4}  using the application outcome as a loss function in the estimation process results in an optimally biased forecasting model that can lead to much better solutions than using standard unbiased ones. In this simple example, because the real-time downward cost is significantly lower than the upward cost, the positive bias is clearly beneficial. However, the complexity of real-world applications prevents us from finding simple solutions as we did before, thereby necessitating more sophisticated approaches like the one proposed in this paper.} 

\section{Convergence Results}\label{sec-conv}

In this section, we discuss some conditions for the convergence of estimators obtained with application-driven joint prediction and optimization. Again, our goal is to obtain the best possible forecast $\hat{y}_t$, but this is completely defined by the parameters $\theta$ since $x_t$ is known. {\color{r2}Let $S_T$ be the set of optimal solutions of \eqref{mod-main-obj}--\eqref{mod-main-rec}, so that  $\theta_T \in S_T$. We will show that any sequence of $\theta_T$, each in the set $S_T$, converges to} the solution set of the actual expected value formulation of the problem (as opposed to the previously presented sampled version). {\color{r4}  We remark that model \eqref{mod-main-obj}--\eqref{mod-main-rec} can be viewed as a particular case of the more general stochastic programming model with equilibrium constraints studied in \cite{shapiro2008stochastic}. In that paper, the authors prove convergence when the data is independent and identically distributed. However, because of the generality of the models in that paper, in order to prove convergence the authors make assumptions (e.g., strong monotonicity) which do not hold in our case. Moreover, as we shall see shortly, we want to extend the result to a non-i.i.d.\ setting. Thus, in our work, we develop a proof that is tailored to our purposes.  }
We will start by describing some assumptions, and then we will state the main theorem. {\color{r4} The proof of the theorem, as well as the proof of the other results in this section, is presented in the Appendix \ref{sec:proof}.}
\begin{assum}
\label{as1}
{\color{r2}There is a unique solution $z_t$ of optimization problem \eqref{mod-main-policy}} for all possible values of $\hat{y}_t$.
\end{assum}
In other words, the problem is always feasible, and the solution {\color{r2}set is a singleton}.
This is not as restrictive as it seems. The feasibility requirement is similar to the classical assumption of complete recourse in stochastic programming. The uniqueness requirement is equivalent to the absence of dual degeneracy in a linear program \citep{borrelli2003geometric}. In this case, the problem in question is dual-degenerate, but it is possible to eliminate this degeneracy by perturbing the objective function---in our case, the vectors $c_p$ and $q_p$---with small numbers that do not depend on the right-hand-side (RHS) of the problem. Thus, the same perturbation is valid for all possible $\hat{y}_t$ \citep{megiddo1989varepsilon}. Another possibility would be resorting to some lexicographic simplex method \citep{nocedal2006numerical}. In this setting, we can define the set-valued function:
\begin{align}
\label{eq:zeta}
\zeta(y) & :=\text{argmin}_{z\in Z}\ G_p(z,y)
\end{align}
From \cite{bohm1975continuity}, we know that if $\zeta(y)$ is a compact set for all $y$ then it is a continuous set-valued function. Moreover, since  $\zeta(y)$ is a singleton for all possible values of $y$, then we treat it as a vector-valued function that is continuous and piecewise affine \citep{borrelli2003geometric}.
\begin{assum}
\label{as_z}
The feasibility set  $Z$ that appears in \eqref{mod-main-policy} is a non-empty and bounded polyhedron.
\end{assum}
Assumption \ref{as_z} is reasonable since this is the set of implementable solutions of the decision-maker, typically representing physical quantities.
\begin{assum}
\label{as_dualQ}
The feasibility set of the dual of the problem that defines $Q_a(z,y)$ in \eqref{mod-main-rec} is non-empty and bounded.
\end{assum}
Note that this set does not depend on $z$ and $y$ since they appear in the RHS of the primal problem. Again, this assumption is akin to a relatively complete recourse assumption applied to the problem defining the outer-level function. 

We state now our main convergence result.
\begin{thm}
\label{thm}
Consider the process given by \eqref{mod-main-obj}--\eqref{mod-main-rec} and {\color{r2}any possible output $\theta_{T} \in S_T$, for each $T$}. Suppose that (i) Assumptions \ref{as1}, \ref{as_z} and \ref{as_dualQ} hold, (ii) the forecasting function $\Psi(\cdot,\cdot)$ is continuous
in both arguments, (iii)  the data process $(X_{1},Y_{1}),\ldots,(X_{T},Y_{T})$
is independent and identically distributed (with $(X,Y)$ denoting a generic element), (iv) the random variable $Y$ is integrable, and (v) the set $\Theta$ is compact and non-empty. Then, {\color{r2}with probability 1,}
\begin{align}
\lim_{T\to\infty}d(\theta_{T},S^{*}) & =0, \label{eq-conv-set}
\end{align}
where $d$ is the Euclidean distance from a point to a set and $S^{*}$ is defined as
\begin{align}
S^{*} & =\textup{argmin}_{\theta\in\Theta}\,\Ea{G_a\big(\zeta(\Psi(\theta,X)),Y\big)}, \label{eq-conv-obj}
\end{align}
with $\zeta(\cdot)$ defined in \eqref{eq:zeta}.
\end{thm}

\begin{rem}
Assumption \ref{as_dualQ} can be replaced by assuming a compact support of $Y$; in this case, $G_a(z, y)$ is a continuous function, where both arguments are defined on compact sets; hence it attains a maximum and is trivially integrable.
\end{rem}

{
\color{r2}
While Theorem \ref{thm}  provides the desired convergence result, condition (iii) of the theorem clearly precludes modeling the situation where the features $x_t$ include (functions of) previous observations $y_{t-1},\ldots,y_{t-k}$. We now extend that result to the case where the features $x_t$ include only lagged observations of $\{y_t\}$.  To do so, suppose that the data process generating $\{Y_t\}_{t=1}^{\infty}$ is a  \textit{stationary ergodic} time series. Stationarity means that the joint distribution of $(Y_1,Y_2,\ldots,Y_k)$ is the same as the joint distribution of $(Y_{t+1},Y_{t+2},\ldots,Y_{t+k})$ for all positive integers $t$ and $k$. It is a standard assumption in the analysis of time series; see, e.g.,  \cite{brockwelldavis:2009} (note that \citealt{brockwelldavis:2009} actually call this notion \textit{strict stationarity}, but elsewhere in the literature it is called just stationarity; see, e.g., \citealt{white:14} or \citealt{billingsley1986}). On the other hand, an ergodic time series is, roughly speaking, one that exhibits a form of ``average asymptotic independence''; a precise definition can be found, for instance, in \cite{white:14}. Statistical tests for stationarity and ergodicity of Markovian processes have been developed by \cite{domowitz1993consistent}. 

\medskip
We can now state a more general version of  Theorem~\ref{thm}:

\begin{thm}
\label{thm2}
Theorem~\ref{thm} is still valid if the assumption that the data process $(X_{1},Y_{1}),\ldots,(X_{T},Y_{T})$
is independent and identically distributed is replaced with the following assumption: $X_t$ is defined as a (measurable) function of $Y_1,\ldots,Y_{t-1}$, and the data process generating $\{Y_t\}_{t=1}^{\infty}$ is a  \textit{stationary ergodic} time series.
\end{thm}

\begin{cor}
\label{cor}
Theorems \ref{thm} and \ref{thm2} imply that if $G_a$ is a cost assessment function describing the ultimate goal of a given practitioner---e.g., the expected cost incurred when using the forecast function $\Psi(\theta,X)$ within a given application---then, under the conditions of Theorem~\ref{thm} or Theorem~\ref{thm2},  any convergent subsequence of the process $\{\theta_T\}_{T=1}^{\infty}$ generated by the estimation process \eqref{mod-main-obj}--\eqref{mod-main-rec} converges to the (not necessarily unique) best forecast model, $\Psi(\theta^*,X)$, in terms of the related application. That is, $\Ea{G_a\big(\zeta(\Psi(\theta^*,X)),Y\big)} \le \Ea{G_a\big(\zeta(\Psi(\theta,X)),Y\big)} \; \forall \theta \in \Theta$. This ensures that our model performs asymptotically better than the classical open-loop approaches. 

\end{cor}

Corollary \ref{cor} highlights that application-driven learning is the best one can do for a fixed triplet of assessment, planning and forecasting functions when the ultimate goal is only minimizing the assessment cost. Although it is possible that other methods lead to the same optimal objective cost, they cannot be better. For other goals, such as minimizing the squared error of forecasts, other methods, such as least squares, will be better since they are inherently aligned with such other goals.
}

\section{Solution Methodology}\label{sec-solution}

In this section, we describe solution methods to estimate the forecasting model within the proposed application-driven closed-loop framework described in \eqref{mod-main-obj}--\eqref{mod-main-rec}. First, we present an exact method based on an equivalent single-level mixed integer linear programming (MILP) reformulation of the bilevel optimization problem \eqref{mod-main-obj}--\eqref{mod-main-rec}. This method uses MILP-based linearization techniques to address the Karush Kuhn Tucker (KKT) optimality conditions of the second level and thereby guarantee the global optimality of the solution in exchange for limited scalability. In the sequence, we describe how to use zero-order methods \citep{conn2009introduction} that do not require gradients to develop an efficient and scalable heuristic method to achieve high-quality solutions to larger instances. These methods will leverage existing optimization solvers, their current implementations and features.

\subsection{MILP-based exact method}\label{sec-bilevel}

Our first approach consists of solving the bilevel problem \eqref{mod-main-obj}--\eqref{mod-main-rec} with standard techniques based on the KKT conditions of the second-level problem \citep{fortuny1981representation}. Thus, the resulting single-level nonlinear equivalent formulation can be reformulated as a MILP and solved by standard commercial solvers. The conversion between the KKT form to the MILP form can be done by numerous techniques \citep{siddiqui2013sos1,pereira2005strategic,fortuny1981representation}, all of which have pros and cons. These techniques are implemented by the open-source package BilevelJuMP.jl  \citep{garcia2023bileveljump} that allows users to formulate and solve bilevel problems in JuMP \citep{lubin2023jump}.

For the sake of completeness, we write the single-level nonlinear reformulation of the bilevel problem \eqref{mod-main-obj}--\eqref{mod-main-rec} in \eqref{mod-kkt-obj}--\eqref{mod-kkt-comp}. For simplicity, in this model, we assume that $Z = \{z| Az \geq h\}$ and that $\Theta$ is polyhedral.
\begin{align}
    \min \limits_{\theta \in \Theta, \hat{y}_t, z_t^{*}, u_t, \pi_t} \quad&  \frac{1}{T}\sum_{t \in \mathbb{T}} \big[ c_a^{\top} z_t^{*} + Q_a(z_t^{*}, {y}_t) \big] \label{mod-kkt-obj} \\
    s.t. \quad& \forall t \in \mathbb{T}: \notag \\
    &\hat{y}_t = \Psi (\theta, x_t) \label{mod-kkt-forecast}\\
        & {W_p}y_t + {H_p}z_t^{*} \geq {b_p} + {F_p} \hat{y}_t \ ; \ \ Az_t^{*}  \geq h \label{mod-kkt-primal} \\
        & {W_p}^{\top}\pi_t = {q_p} \ ; \ \ H_p^{\top}\pi_t + A^{\top}\mu_t = {c_p} \ ; \ \ \pi_t,\mu_t \geq 0  \label{mod-kkt-dual} \\
        & \pi_t \perp {W_p}u_t  + {H_p}z_t^{*} - {b_p} - {F_p} \hat{y}_t \ ; \ \ \mu_t \perp {A}z_t^{*} - h \label{mod-kkt-comp}
\end{align}
Equations \eqref{mod-kkt-obj} and \eqref{mod-kkt-forecast} are the same as \eqref{mod-main-obj} and \eqref{mod-main-forecast}. %This model uses the same definition of (\ref{mod-main-ope}). 
Equation \eqref{mod-main-policy} was replaced by \eqref{mod-kkt-primal}--\eqref{mod-kkt-comp}. Constraints \eqref{mod-kkt-primal}, \eqref{mod-kkt-dual}, and \eqref{mod-kkt-comp} represent, respectively, the primal feasibility constraints, the dual feasibility constraints, and the complementarity constraints.

\subsection{Scalable heuristic method}\label{sec-zero}

The proposed class of methods will make extensive use of the way of thinking described in Figure \ref{fig-closedloop}. In other words, the core algorithm decomposes the problem as detailed in Algorithm \ref{alg1}.

\begingroup
\renewcommand{\baselinestretch}{0.8}
\begin{algorithm}[h]
\SetAlgoLined
\KwResult{Optimized $\theta$}
 Initialize $\theta$ \;
 \While{Not converged}{
  Update $\theta$\;
  \For{$t \in \mathbb{T}$}{
    Forecast: $\hat{y}_t \leftarrow \Psi (\theta, x_t)$\;
    Plan Policy: $z^{*}_t \leftarrow \argmin_{z \in Z } G_p(z, \hat{y}_t)$\;
    Cost Assessment: ${cost}_t \leftarrow  G_a(z^{*}_t, {y}_t)$
  }
  Compute cost: $cost(\theta) \leftarrow  \sum_{t \in \mathbb{T}}({cost}_t)$
 }
\caption{Pseudo-algorithm}
\label{alg1}
\end{algorithm}
\endgroup
We call this method a pseudo-algorithm because a few steps are not fully specified \textit{a priori}, namely \textit{Initialization}, \textit{Update}, and \textit{Convergence} check, allowing for a wide range of possible specifications.
\textit{Initialization} can be as simple as $\theta$ receiving a vector of zeros, which might not be good if the actual algorithm is a local search. One alternative that will be applied in the case study section is the usage of traditional models as starting points, for instance, the ordinary least squares.
In the case study we will initialize the algorithm with the LS estimate. This guarantees that the algorithm will return at most the same cost as the open-loop framework in the training sample.
There are many possibilities for the \textit{convergence} test. For instance, iteration limit, time limit, the variation of the objective function value, and other algorithm-specific tests. Finally, the \textit{update} step depends on the selected concrete algorithm that ultimately minimizes the non-trivial $cost(\theta)$ function.

We will focus on a derivative-free local search algorithm named Nelder-Mead \citep{conn2009introduction}.
Notwithstanding, it is relevant to highlight the generality of the proposed pseudo-algorithm. For instance, gradient-based algorithms could also be developed based on numerical differentiation and automatic differentiation \citep{nocedal2006numerical}.
In this context, gradient calculation would enable the usage of Gradient Descent and BFGS-like algorithms \citep{nocedal2006numerical}.

The main features of the above-proposed pseudo-algorithm are: 1) it is suitable for parallel computing (the loop in the sample $\mathbb{T}$ is intrinsically decoupled); 2) each step is based on a deterministic LP defining the second-level variables in \eqref{mod-main-policy}, suitable for off-the-shelf commercial solvers that find globally optimal solutions in polynomial time; 3) each inner step can significantly benefit from warm-start processes developed in linear programming solvers (e.g., the dual simplex warm-start is extremely powerful, and many times only a handful of iterations will be needed in comparison to possibly thousands of iterations if there were no warm-start, cf. \cite{nocedal2006numerical}). It is worth emphasizing that the aforementioned feature 2) allows for an exact (always optimal) description of the second-level problem. 
In our approach, we keep the second level exact and face the challenge of optimizing a nonlinear problem on the upper level. In contrast, \cite{munoz2022bilevel} choose to relax the complementarity constraints and deal with a nonlinear program lacking the benefits of the above-mentioned features 1) to 3). Finally, note that this heuristic approach allows for a wider set of forecast models, such as Neural Networks and other machine learning models, as it only requires that a pointwise forecast can be obtained and its performance evaluated by the cost function for a given trial solution (parameters).

One caveat is that variations on $\theta$ can lead to possibly infeasible results for the \textit{Policy Planning} and \textit{Cost Assessment} optimization problems. Consequently, we require complete recourse for such problems. In cases where this property does not hold, it is always possible to add artificial (slack) variables with high penalty costs in the objective function to keep the problem feasible. 

\section{Application-Driven Load Forecasting and Reserve Sizing} \label{sec-power}

In this work, we focus on the energy and reserve scheduling problem of power systems \citep{Chen2014, kirschen2018}. 
We apply the framework of Section \ref{sec-theory} to obtain the best joint conditional point-forecast for the vector of nodal demands, $\hat{D}_t$, and vectors of up and down zonal or nodal reserve requirements, $\hat{R}_t^{(up)}$ and $\hat{R}_t^{(dn)}$. Hence, $\hat{y}_t = (\hat{D}_t, \hat{R}_t^{(up)}, \hat{R}_t^{(dn)})$. The forecast vector of nodal demands represents, e.g., the next hour operating point target that system operators and agents should comply with. The up- and down-reserve requirements represent a forecast of the system's resource availability (or security margins), defined per zone or node, allowing the system to withstand load deviations. Note that we can think of loads as a general net load that corresponds to load minus non-dispatchable (e.g., variable renewable) generation. {\color{r4}The forecasts are obtained from parametric functions that depend on explanatory variables, $x_t$, analogous to equation \eqref{mod-main-forecast}:
\begin{align}
& \hspace{52px} \hat{D}_t = \Psi_D (\theta_D, x_t) \label{mod-pow-fore-d}\\
& \hspace{52px}   \hat{R}_t^{(up)} = \Psi_{R^{(up)}} (\theta_{R^{(up)}}, x_t) \label{mod-pow-fore-up}\\
& \hspace{52px}   \hat{R}_t^{(dn)} = \Psi_{R^{(dn)}} (\theta_{R^{(dn)}}, x_t) \label{mod-pow-fore-dn}
\end{align}
Hence, $\theta = (\theta_D, \theta_{R^{(up)}}, \theta_{R^{(dn)}})$.

Now we describe two main types of input data.
First, we have historical data, $\{y_t, x_t\}_{t\in \mathbb{T}}$, including dependent variables ${y}_t = ({D}_t, {R}_t^{(up)}, {R}_t^{(dn)})$ (only ${D}_t$ will be in the model) and explanatory variables, for instance, lags of demand, ${x}_t = (D_{t-1},...,D_{t-k})$, but also other exogenous covariates such as climate and weather indices (or, e.g., external forecasts). Second, we consider system physical and market data, in vector form, as follows.
Generating units data (vectors): maximum generation capacities $K$, dispatch costs or offers $c$, maximum up- and down-reserve capacities $\bar{r}^{(up)}$ and  $\bar{r}^{(dn)}$, up- and down-reserves costs $p^{(up)}$ and $p^{(dn)}$; network data: vector of transmission line capacities $F$, network sensitivity matrix $B$, describing the network topology and physical laws of electric circuits; load-shed and spillage penalty costs $\lambda^{LS}$ and $\lambda^{SP}$.
For a simple matrix representation of the problem, we define $e$ to be a vector with one in all entries and an appropriate dimension.
We also need to define matrices to describe the bus or zone where each resource is located.
$M$ is an incidence matrix with buses in rows and generators in columns that is one when the generator lies in that bus and zero otherwise. Similarly, $N$ is an incidence matrix with generators in columns and reserve zones in rows, which is one if the generator lies in that zone and zero otherwise.}

The input data describing the system characteristics can be provided under two perspectives: 1) under the perspective of the actual \emph{ex-post} (or assessed/implemented) operation, i.e., based on the most accurate system's description for {\color{r2}optimizing} the function $G_a$ defined in \eqref{mod-main-obj}; and 2) under the \emph{ex-ante} planning perspective, $G_p$, which is accounted for in \eqref{mod-main-policy} based on system operator's description of the system considered in the dispatch optimization tool. While the former has already been listed in the previous paragraph, the latter uses the same symbols but with a tilde above, i.e., $\tilde{c},\tilde{p}^{(up)}, \tilde{p}^{(dn)}, \tilde{K}, \tilde{B}$, etc.

{\color{r4}
Given the system data and the forecasts $\hat{y}_t = (\hat{D}_t, \hat{R}_t^{(up)}, \hat{R}_t^{(dn)})$ we can obtain the plan $z_t^*$ by optimizing the planning scheme, $G_p$. The optimal plan, $z_t^*$, is comprised of the generation setting point, $g^*_t$, and allocated up and down reserves, $r_t^{(up)*}, r_t^{(dn)*}$, that is, $z^* = (g^*, r_t^{(up)*}, r_t^{(dn)*})$. The optimization of the planning phase defining function, $G_p$ from \eqref{mod-main-policy}, can be written as follows:
\begin{align}
& \hspace{42px} \Big( g_t^* , r_t^{(up)*}, r_t^{(dn)*}\Big) \in  \argmin_{\substack{\hat{g_t}, \hat{\delta}_t^{LS}, \hat{\delta}_t^{SP},\\\hat{r_t}^{(up)}, \hat{r_t}^{(dn)}}}  \big[\tilde{c}^{\top} \hat{g_t} + \tilde{p}^{(up)\top} \hat{r_t}^{(up)} + \tilde{p}^{(dn)\top} \hat{r_t}^{(dn)} +  
            \tilde{\lambda}^{LS}\hat{\delta}_t^{LS} + \tilde{\lambda}^{SP}\hat{\delta}_t^{SP} \big]  \label{mod-pow-pl-obj}\\
& \hspace{155px} s.t. \hspace{30px} e^\top (M \hat{g_t} - \hat{\delta}_t^{SP}) = e^\top (\hat{D}_t - \hat{\delta}_t^{LS}) \label{mod-pow-pl-bal}\\
& \hspace{204px} -\tilde{F} \leq \tilde{B}(M\hat{g_t} + \hat{\delta}_t^{LS} - \hat{D}_t - \hat{\delta}_t^{SP}) \leq \tilde{F} \label{mod-pow-pl-net}\\
& \hspace{204px}  N \hat{r_t}^{(up)} = \hat{R}_t^{(up)} \label{mod-pow-pl-balup}\\
& \hspace{204px} N \hat{r_t}^{(dn)} = \hat{R}_t^{(dn)} \label{mod-pow-pl-baldn}\\
& \hspace{204px} \hat{g_t} + \hat{r_t}^{(up)} \leq \tilde{K} \label{mod-pow-pl-ub}\\
& \hspace{204px} \hat{g_t} - \hat{r_t}^{(dn)} \geq 0 \label{mod-pow-pl-lb}\\
& \hspace{204px} \hat{r_t}^{(up)} \leq {\bar{r}}^{(up)} \label{mod-pow-pl-upub}\\
& \hspace{204px} \hat{r_t}^{(dn)} \leq {\bar{r}}^{(dn)} \label{mod-pow-pl-dnub}\\
& \hspace{204px} \hat{g_t}, \hat{r_t}^{(up)}, \hat{r_t}^{(dn)}, \hat{\delta}_t^{LS}, \hat{\delta}_t^{SP} \geq 0. \label{mod-pow-pl-zero}
\end{align}
Although there are some variations around the world, operators solve deterministic problems similar to \eqref{mod-pow-pl-obj}--\eqref{mod-pow-pl-zero} in their planning phase at each period to define a generation and reserve schedule for the next period.
The objective function, \eqref{mod-pow-pl-obj}, aims to minimize a weighted sum of the generation setting point, $\hat{g_t}$, allocated up and down reserves, $r_t^{(up)*}, r_t^{(dn)*}$, load shedding, $\hat{\delta}_t^{LS}$, and generation spillage (excess), $\hat{\delta}_t^{SP}$. Constraint \eqref{mod-pow-pl-bal} ensures that total generation minus the excess will match the total forecast demand minus the load-shed. Constraint \eqref{mod-pow-pl-net} enforces that energy flows in the transmission lines will be within their limits. Expressions \eqref{mod-pow-pl-balup} and \eqref{mod-pow-pl-baldn} ensure that the forecast reserve requirements ($\hat{R}^{(up)}_t$ and $\hat{R}^{(dn)}_t$), which are considered as input data for this model, must be allocated among generators in the form of up and down reserves ($\hat{r_t}^{(up)}$ and $\hat{r_t}^{(dn)}$). Constraints \eqref{mod-pow-pl-ub} and \eqref{mod-pow-pl-lb} limit the scheduled generation and reserves range (up and down) to generators' physical generation limits. Constraints \eqref{mod-pow-pl-upub} and \eqref{mod-pow-pl-dnub} limit the maximum amount of reserves that can be allocated in each generating unit, and \eqref{mod-pow-pl-zero} ensures positiveness of the generation, up- and down-reserves, load-shed and generation spillage decision vectors.

Finally, given the planned decision (generation setting points, $g_t^*$, and reserve schedules, $\hat{r_t}^{(up)*}, \hat{r_t}^{(dn)*}$), and the realization of the demand, $D_t$, we can assess the actual system cost. This final evaluation is done by the optimization problem defined by $G_a$ in \eqref{mod-main-obj}. The problem corresponds to the phase in which system operators adjust their planned decisions to cope with real-time events. Thus, we define the evaluation model as:
\begin{align}
& \min_{\substack{g_t, \delta_t^{LS}, \delta_t^{SP}}} 
        \frac{1}{T}\sum_{t \in \mathbb{T}}
            \big[c^{\top} g_t^* + p^{(up)\top} \hat{r_t}^{(up)*} + p^{(dn)\top} \hat{r_t}^{(dn)*} + 
            \lambda^{LS}\delta_t^{LS} + \lambda^{SP}\delta_t^{SP} \big]  \label{mod-pow-op-obj}\\
& \hspace{15px} s.t. \hspace{18px} \forall t \in \mathbb{T}: \notag \\
& \hspace{52px} e^{\top} (M g_t - \delta_t^{SP}) = e^\top (D_t - \delta_t^{LS}) \label{mod-pow-op-bal}\\
& \hspace{52px} -F \leq B(M g_t + \delta_t^{LS} - D_t - \delta_t^{SP}) \leq F  \label{mod-pow-op-net}\\
& \hspace{52px} g_t^* - r_t^{(dn)*}\leq g_t \leq g_t^* + r_t^{(up)*}  \label{mod-pow-op-adj}\\
& \hspace{52px} \delta_t^{LS}, \delta_t^{SP}, \hat{R}_t^{(up)}, \hat{R}_t^{(dn)}, g_t \geq 0 \label{mod-pow-op-zeeo}
\end{align}
where, the objective function, \eqref{mod-pow-op-obj}, considers fixed costs of planned decisions, which are input data, and minimizes the costs of imbalances (load-shed, $\delta_t^{LS}$, and generation spill, $\delta_t^{SP}$). Constraints \eqref{mod-pow-op-bal} and \eqref{mod-pow-op-net} are similar to \eqref{mod-pow-pl-bal} and \eqref{mod-pow-pl-net}. Constraint \eqref{mod-pow-op-adj} limits the \emph{ex-post} generation to respect the operation range defined by the \emph{ex-ante} planned decision. Constraint (\ref{mod-pow-op-zeeo}) ensures the positiveness of imbalance ($\delta_t^{SP}, \delta_t^{LS}$), generation, and reserve requirement variables.

Combining 
\eqref{mod-pow-fore-d}--\eqref{mod-pow-fore-dn},
\eqref{mod-pow-pl-obj}--\eqref{mod-pow-pl-zero}
and
\eqref{mod-pow-op-obj}--\eqref{mod-pow-op-zeeo} we completely specialize the closed-loop application-driven learning framework proposed in \eqref{mod-main-obj}--\eqref{mod-main-rec}. The learning framework results in a joint load and reserve requirements forecasting model that results in the minimum operational cost under observed data. Based on the proposed framework, optimally biased forecasts for system loads and reserve requirements ($\hat{y}_t = \Psi(\theta, x_t)$) can be used to guide deterministic planning models towards more economical energy and reserve scheduling decisions. In this context, we can say that the forecasts obtained from our application-driven learning framework are not only contextually adapted, as they are functions of all data in $x_t$, but also system-informed, as they account for both the planning and implementation processes ($G_p$ and $G_a$, respectively).
Clearly, $G_a$ and $G_p$ are different functions as they represent different steps in the application pipeline.
We present the full bilevel formulation in the Appendix \ref{app-full}.}

In the case studies, we will compare a few variants of this problem by optimizing all or some of the parameters $\theta_D$, $\theta_{R^{(up)}}$, and $\theta_{R^{(dn)}}$. As described in the contributions section, two novel methods are the one that optimizes all parameters and the one that gets a fixed $\theta_D$ and only optimizes reserves ($\theta_{R^{(up)}}$, and $\theta_{R^{(dn)}}$). 

{\color{r4} 
\section{Linearly biased forecasts: a benchmark model}\label{sec-naive}

As mentioned in \cite{CAISO}, operators have relied on \emph{ad hoc} upward-biased demand forecasts. However, there is no clear reference to how this \textit{ad hoc} bias is obtained. In order to thoroughly test our method (which will automatically bias the forecasts optimally), we briefly present a novel benchmark based on a possible implementation of this industry practice. 

Given a forecast $\hat{y}_t$ obtained by non-application-driven methods (e.g., LS), we define a linearly biased forecast as $\hat{y}_t^{LB} = \alpha^* \hat{y}_t$ where 
\begin{align}
    \alpha^* \in \argmin \limits_{\alpha \in \mathbb{R},z_t^{*}} \quad&  \frac{1}{T} \sum_{t \in \mathbb{T}} G_a(z_t^{*}, {y}_t) \label{mod-naive-obj} \\
    s.t. \quad& z_t^{*} \in \argmin_{z \in Z }  G_p(z, \alpha\hat{y}_t)\ \forall t \in \mathbb{T}, \label{mod-naive-policy}
\end{align}\label{model}
This model can be seen as a constrained version of the application-driven framework, where a simple optimized bias scaling factor is applied to the unbiased forecast. As $\alpha$ is a scalar, it can be found through a simple grid search.
As this model was not previously described in the literature, its description and quantitative analyses can also be seen as an interesting practical yet secondary contribution of this work. In Section~\ref{sec-realdata2}, we present numerical results for the implementation of this model in a large-scale system and compare them with those obtained with the application-driven approach we propose in this paper. Hereinafter, this approach is referred to as the linear-bias benchmark approach or linearly biased forecast.
}

\section{Case Studies}\label{sec-simulations}

This section presents case studies to demonstrate the methodology's applicability and how the closed-loop framework can outperform the classic open-loop scheme in multiple variants of the load forecasting and reserve sizing problem defined in Section \ref{sec-power}. {\color{r4} First, we present some considerations about the studied systems and models in Section \ref{sec-system} and \ref{sec-ar}.
Second, in Section \ref{sec-bxh}, we show that the Heuristic method of Section \ref{sec-zero} can achieve close to optimal solutions in a fraction of the time required by the Exact method of Section \ref{sec-bilevel}.
Third, in Section \ref{sec-heuristic}, we study the estimated parameters' and forecasts' empirical properties and contrast them with the classical least squares (LS) estimators.
Finally, in Sections \ref{secVeryLarge} and \ref{sec-realdata2}, we show that the heuristic algorithm finds good quality local-optimal parameters systematically outperforming the LS open-loop benchmark for instances far larger and more realistic than those solved in previously reported works tackling closed-loop bilevel frameworks. Section \ref{sec-realdata2} will contrast the models from Sections \ref{sec-power} and \ref{sec-naive}. An additional study focused on the scalability of the method is presented in the Appendix \ref{sec:scale}.}

\subsection{Power systems cases and datasets}\label{sec-system}

We consider multiple power system cases throughout this section. The first, considered in Sections \ref{sec-bxh} and \ref{sec-heuristic}, is a single bus system defined by us, with 1 zone, 1 load (with a long-term average of $6$) and 4 generators (with capacities $5, 5, 2.5,2.5$ and costs $1, 2, 4, 8$).
Then, in Sections \ref{secVeryLarge} and \ref{sec-realdata2}, we consider realistic power systems with more than 6000 buses from PG-LIB-OPF \citep{babaeinejadsarookolaee2019power}.

Henceforth, we will refer to the test systems by their number of buses. We only considered as loads the buses with positive demand in the original files.
Deficit and generation curtailment costs were defined, respectively, as $8$ and $3$ times the most expensive generator cost. All generators were allowed to have up to $30\%$ their capacity allocated to reserves, and their reserve allocation costs were set to $30\%$ their nominal costs. We only considered the linear component of the generators' costs in all instances. 

In all datasets, except in Section
\ref{sec-realdata2}, we used demand values as the long-term average of AR(1) (autoregressive of order 1) processes for each bus. The AR(1) coefficients were set to $0.9$, and the AR(0) coefficients were set so that we get the desired long-term averages. For the sake of simplicity, load profiles were generated independently. The coefficient of variation of all simulated load stochastic processes was $0.4$. Negative demands were truncated to zero, although they could represent an excess renewable generation.
In Section \ref{sec-realdata2}, we demonstrate the performance of our approach when considering real demand data.

\subsection{Studied models and notation}\label{sec-ar}

In the following sections, we will consider a linear autoregressive forecast model. This allows us to compare the results with the true data model in a controlled environment. Therefore, unless otherwise mentioned, the model used to forecast loads in each node is the following:
\begin{align}
\hat{D}_t = \Psi_D (\theta_D, x_t) = \theta_D(0) + \theta_D(1) D_{t-1}, \label{mod-ar1-dem}
\end{align}
For the single bus case, we set the ``real", or population, values as $\theta_D(0) = 0.6$ and $ \theta_D(1) = 0.9$, resulting in the long-term average $\theta_D(0) / (1- \theta_D(1)) = 6$, defined in section \ref{sec-system}. {\color{r2}Note that such a choice of coefficients ensures that the  input process is stationary and ergodic (see, e.g., \citealt[p.495]{billingsley1986}), thereby satisfying the condition in Theorem~\ref{thm2}.} 
Because the stochastic model for loads is homoscedastic, we set the reserve models to  $AR(0)$ -- i.e., a number that does not depend on previous values of reserves -- since it is customary to set the reserves just in terms of variability of loads:
\begin{align}
\hat{R}_t^{(up)} = \Psi_{R^{(up)}} (\theta_{R^{(up)}}, x_t) = \theta_{R^{(up)}}(0), \label{mod-ar10-rup} \\
\hat{R}_t^{(dn)} = \Psi_{R^{(dn)}} (\theta_{R^{(dn)}}, x_t) = \theta_{R^{(dn)}}(0). \label{mod-ar10-rdn}
\end{align}

\subsection{Exact vs heuristic method comparison}\label{sec-bxh}

In this first experiment, we aim to compare the exact and heuristic methods to check the quality of the latter. We consider the single-bus test system so that the exact method can reach global optimal solutions. We use a Dell Notebook (Intel i7 8th Gen with 4 cores at 1.99Ghz, 16Gb RAM).

We solved $10$ instances for each $T \in \{15, 25, 50, 75\}$.
All instances solved with the exact method converged within a gap lower than $0.1\%$ using the Gurobi solver
\citep{gurobi}
or stopped after two hours.
The heuristic method was terminated when the objective function presented a decrease lower than $10^{-7}$ between consecutive iterations. We used a Nelder-Mead implementation found in \cite{mogensen2018optim}.
To compare the results, we plotted the ratio of objective values in Figure \ref{fig-scatt}(a) and the time ratio in 
Figure \ref{fig-scatt}(b).
We can observe that the heuristic method achieves high-quality solutions for almost all instances. Although the exact method is competitive for $T \in \{15, 25\}$, the heuristic method is much faster with an average solve time of $4.4s$, for $T = 50$, and $5.9s$, for $T = 75$, compared to $1200s$ and $6670s$ for the exact method. The exact method did not converge for four instances with $T = 75$ after two hours.
\begin{figure}[h!]
\begin{center}
\subfloat[In-Sample objective value comparison]{
\includegraphics[height=2.0in]{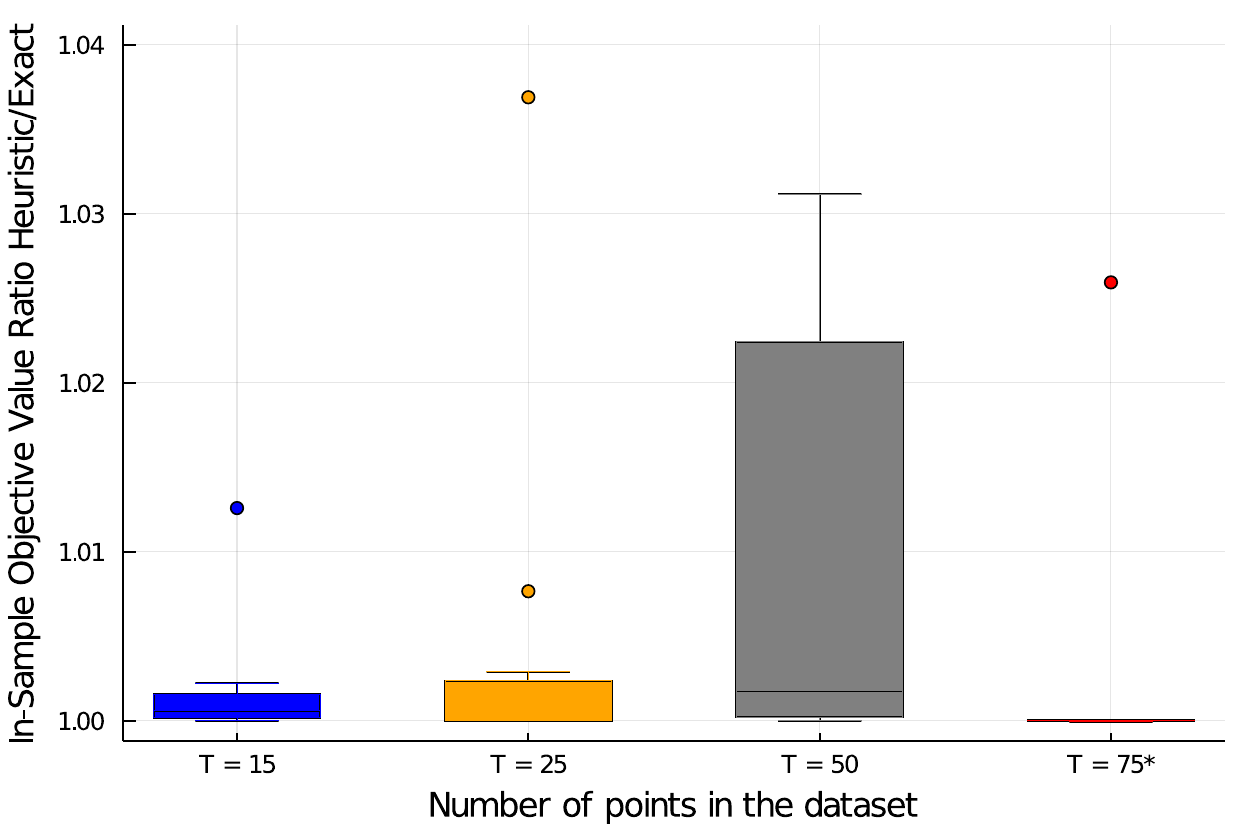}
}
\subfloat[Running time comparison]{
\includegraphics[height=2.0in]{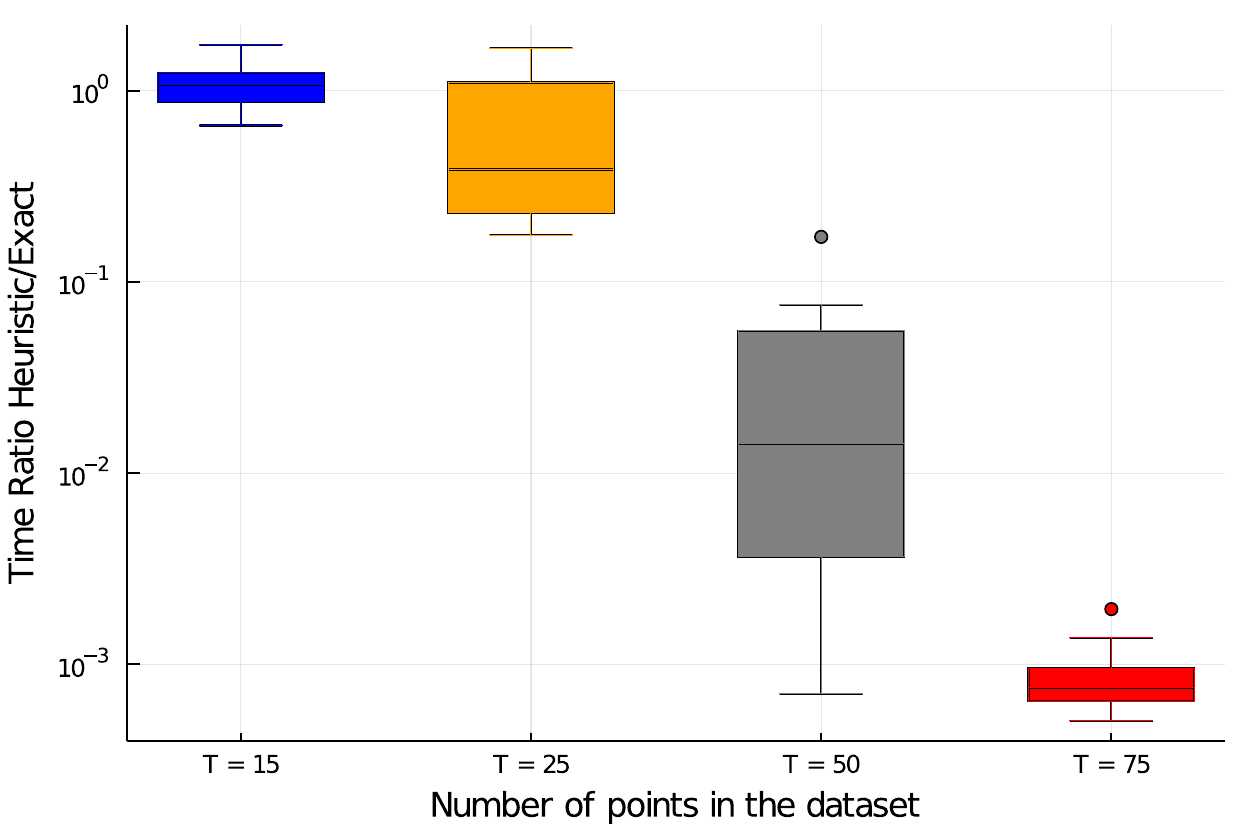}
}
\end{center}
\caption{(a) Objective of Heuristic method divided by the objective of Exact method for the same datasets. *Four problems not considered for $T = 75$: the exact method found no solution in the given time.\\ (b) Time to solve the same problems (in log scale): Heuristic method divided by Exact method.} \label{fig-scatt}
\end{figure}

\subsection{Asymptotic behavior and biased estimation}\label{sec-heuristic}

Now we focus only on the heuristic method to analyze how the estimates behave with respect to the variation in dataset size. We will see that they actually converge in our experiments. Moreover, we empirically show through out-of-sample studies that using the closed-loop model is strictly better than the open-loop one, provided we have a reasonable dataset size in the training step. It will be possible to see that a method with too many parameters might overfit the model for reduced dataset sizes and not generalize well enough. We use a Dell Notebook (Intel i7 8th Gen with 4 cores at 1.99Ghz, 16Gb RAM).
From now on, we will use the following nomenclature and color code to refer to the different models:
\begin{itemize}
	\item LS-Ex (red): This is the benchmark model representing the classical open-loop approach. It uses LS to estimate demand and an exogenous reserve requirement.
	\item LS-Opt (blue): This is a partially optimized model, where least squares are used to estimate demand and only reserve requirements are optimized with the application-driven framework.
    \item Opt-Ex (yellow): This is also a partially optimized model, where demand is optimized, whereas reserve requirements are still exogenously defined. This model is not particularly meaningful in practice. We show it in some studies for completeness.
	\item Opt-Opt (green):  This is the fully optimized model, where both demand forecast and reserve requirements are jointly optimized.
\end{itemize}

\noindent Both (closed-loop) methods LS-Opt and Opt-Opt are novel contributions proposed in this work, to be contrasted with the benchmark (open-loop) method, LS-Ex. For didactic purposes, in all cases tested in this section, up and down reserve requirements were defined as +/- $1.96$ standard deviations, respectively, of the estimated residuals from the LS demand forecast. Note, however, that other exogenous \emph{ad hoc} rule could be used \citep{NRELguidelines}. 

We empirically compare and analyze the convergence of the four demand and reserve requirement forecast models mentioned above. We varied the dataset size used in the estimation process from $50$ to $1000$ observations. For each dataset size, we performed $100$ trial estimations, with different datasets generated from the same process, to study the convergence. To evaluate the out-of-sample performance of each one of the $100$ estimates for each dataset size, we compute the objective function, $G_a$ in \eqref{mod-main-obj}, for a single fixed dataset with $10000$ new observations (generated with the same underlying process but different from all other data used in the estimation/training phase). In the following plots, lines represent mean values among the $100$ estimated costs with the in-sample or out-of-sample data, and shaded areas represent the respective $10\%$ and $90\%$ quantiles.

Figures \ref{fig-testquant} (a) and (b) depict the costs in the out-of-sample data. Hence, they measure how well the models generalize to data it has never seen before. We can see that the models allowing more parameters to be endogenously optimized perform much better than models with exogenously defined forecasts. Thus, we see that the application-driven learning framework works successfully on out-of-sample data when estimated with datasets larger than $150$. However, we note that these steady improvements require more data than the classic exogenous models, as shown in Figure \ref{fig-testquant} (b). Between $50$ and $120$ points, the model with more optimization flexibility, Opt-Opt, exhibits a more significant cost variance. This is due to excessive optimization in a small dataset that led to overfitting and poor generalization. Note that, in this work, we did not consider any regularization procedure to avoid this issue. However, our optimization-based framework is suitable for well-known shrinkage operators \citep{tibshirani2011regression} that can be readily added in the objective function \eqref{mod-main-obj}. 

\begin{figure}[h!]
\subfloat[Dataset size from 100 to 1000.]{
\includegraphics[height=2.0in]{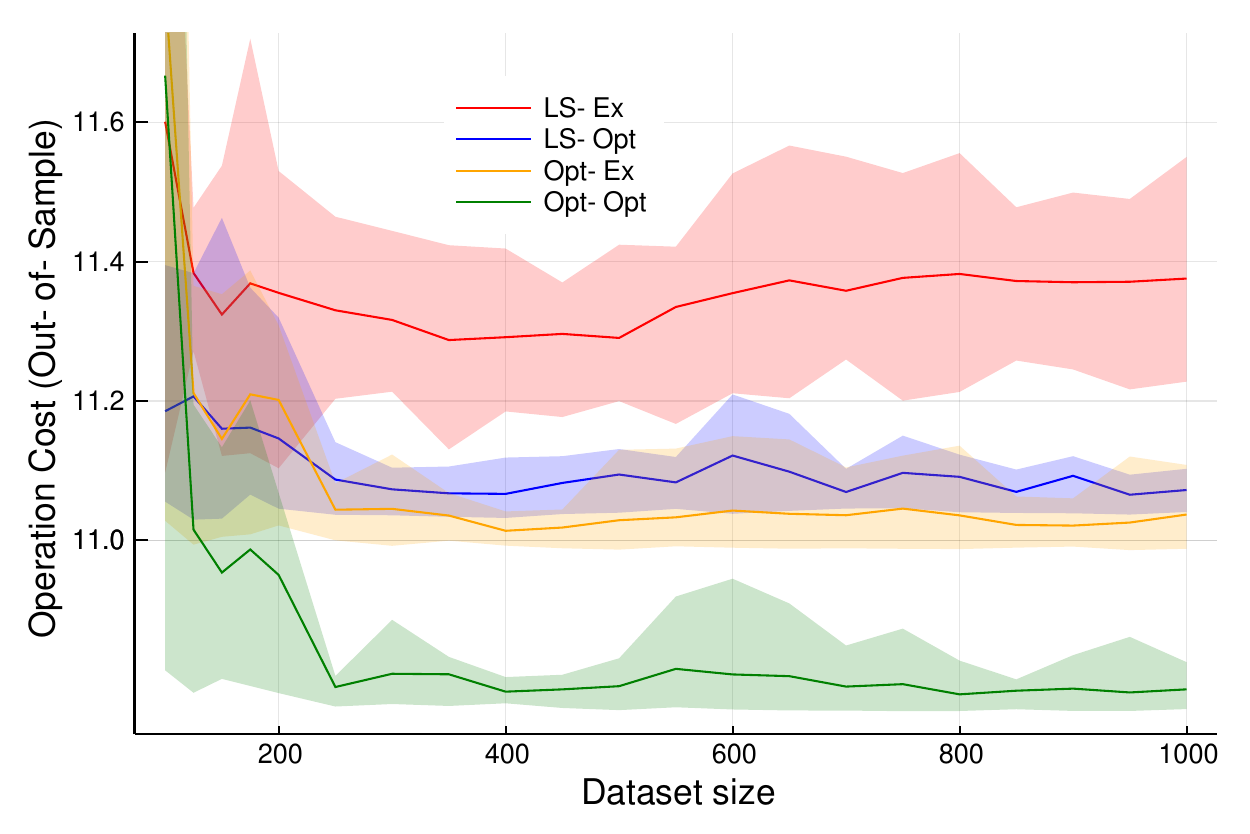}
}
\subfloat[Dataset size from 50 to 200.]{
\includegraphics[height=2.0in]{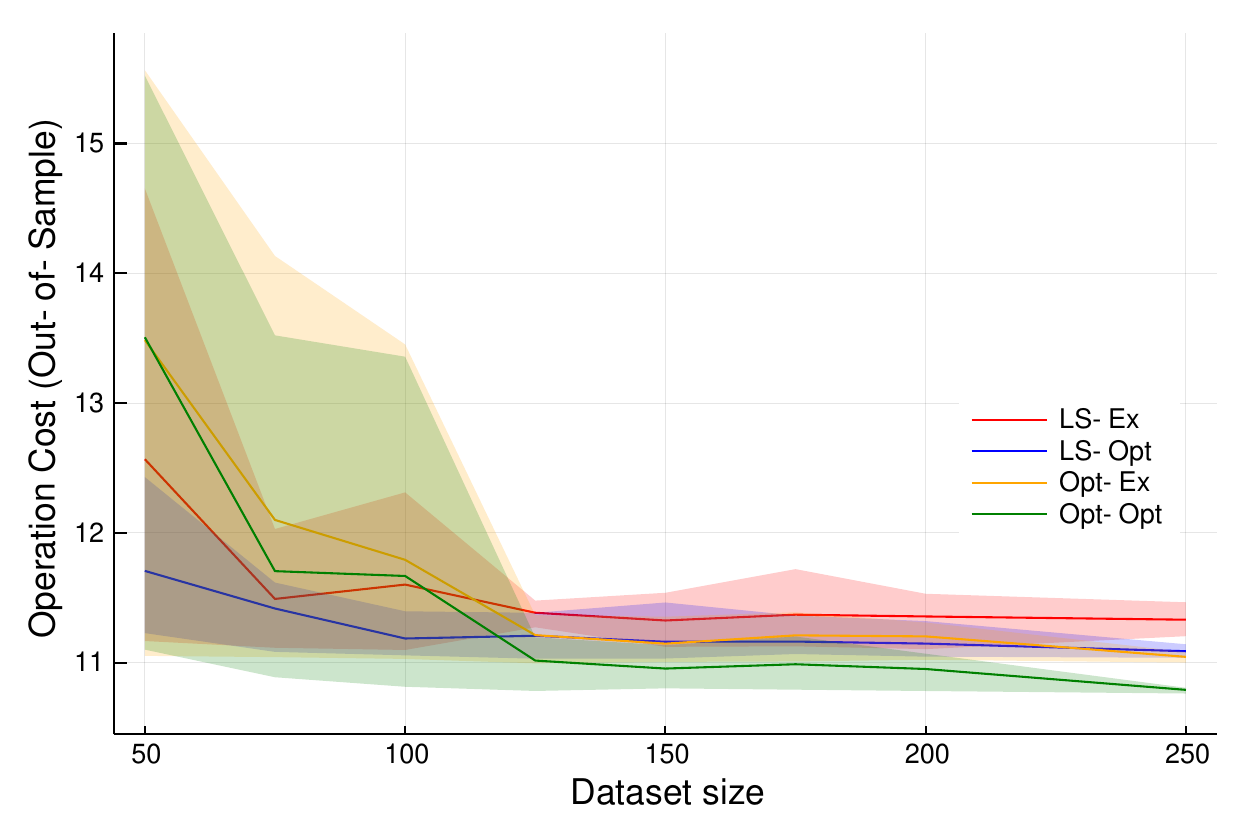}
}
\caption{Out-of-sample average operation cost versus (in-sample) dataset size. Lines represent the average of the 100 estimation trials. Shaded areas represent the $10\%$ and $90\%$ quantiles. All trials are evaluated on a single out-of-sample dataset with size $10,000$ observations.} \label{fig-testquant}
\end{figure}

Figure \ref{fig-demand coef-case1} shows how the estimated parameters behave as functions of the estimation dataset size. In Figure \ref{fig-demand coef-case1} (a) and (b), we can see that the load model parameters are indeed converging to long-run values. It is also clear to see the bias in those parameters. The constant term is greatly increased while the autoregressive coefficient is slightly reduced. Ultimately, this leads to a larger forecast value, which can be interpreted as the application risk adjustment due to the asymmetric imbalance penalization costs (load-shed is much higher than the spillage cost). Thus, the Opt-Opt model will do the best possible to balance these costs, thereby prioritizing the load shed by increasing the forecast level. The fixed reserves model (Opt-Ex) is less biased because the fixed reserves constrain how much the load model can bias due to the risk of not having enough reserves to address lower demand realizations. Note that the red (LS-Ex) is on top of the purple (LS-Opt) since both use the same LS estimates for demand, which exhibits the lowest variance.
\begin{figure}[h!]
\subfloat[$\theta_D(0)$, from \eqref{mod-ar1-dem}]{
\includegraphics[height=2.0in]{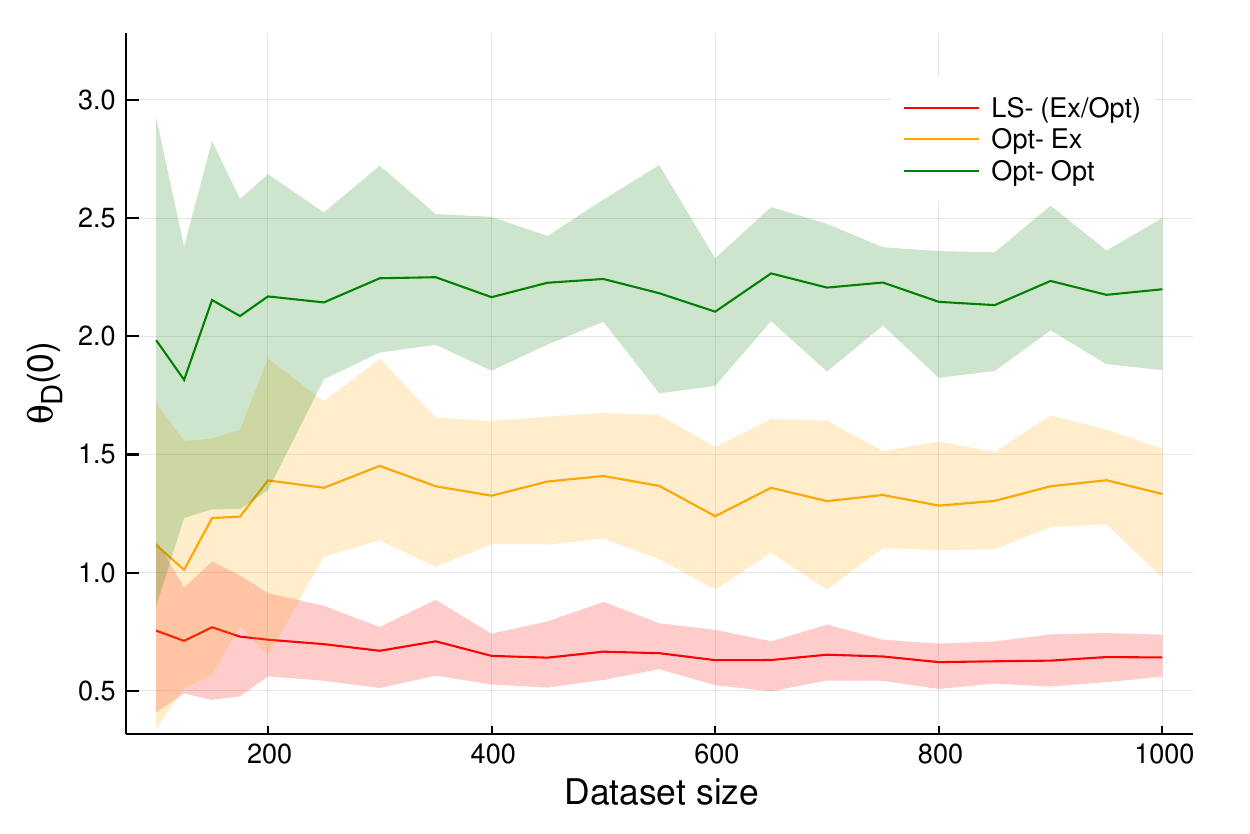}
}
\subfloat[$\theta_D(1)$, from \eqref{mod-ar1-dem}]{
\includegraphics[height=2.0in]{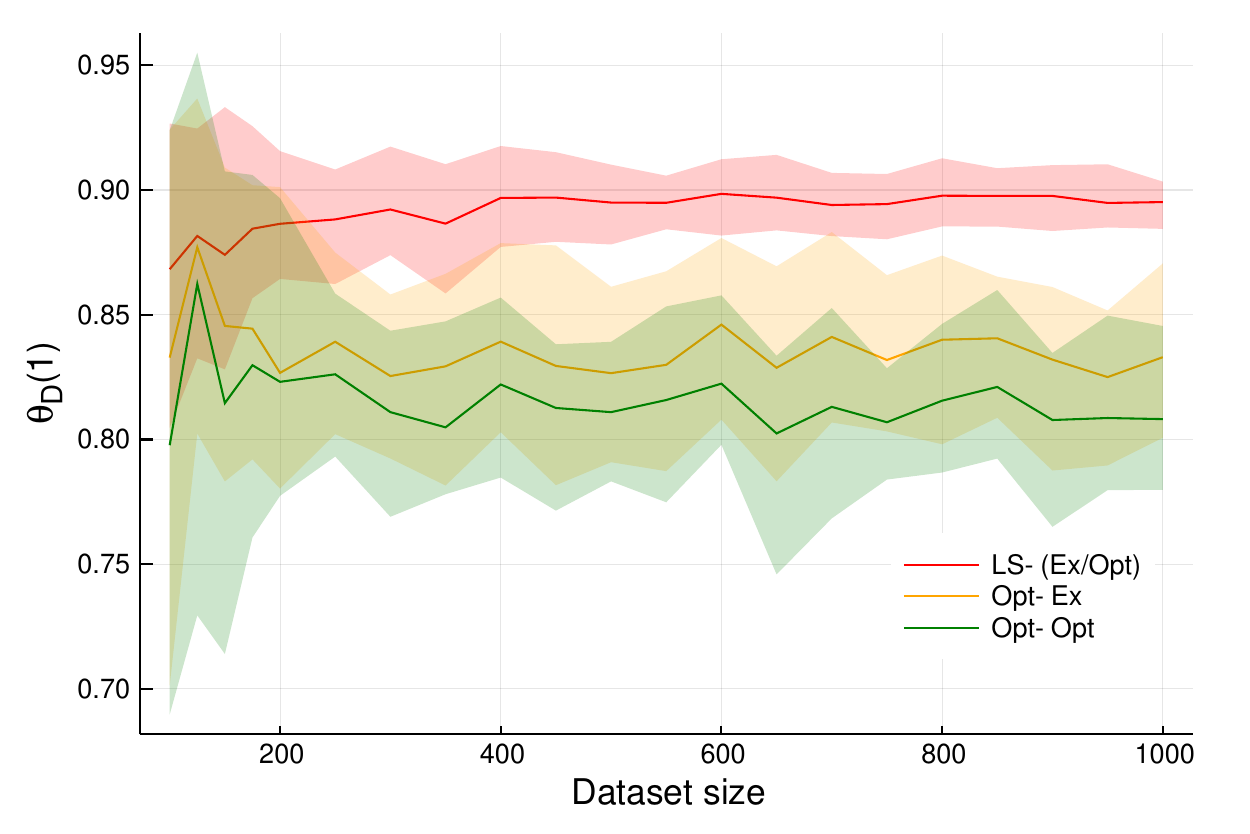}
}
\
\subfloat[Reserve up, $\theta_{R^{(up)}}(0)$, from \eqref{mod-ar10-rup}]{
\includegraphics[height=2.0in]{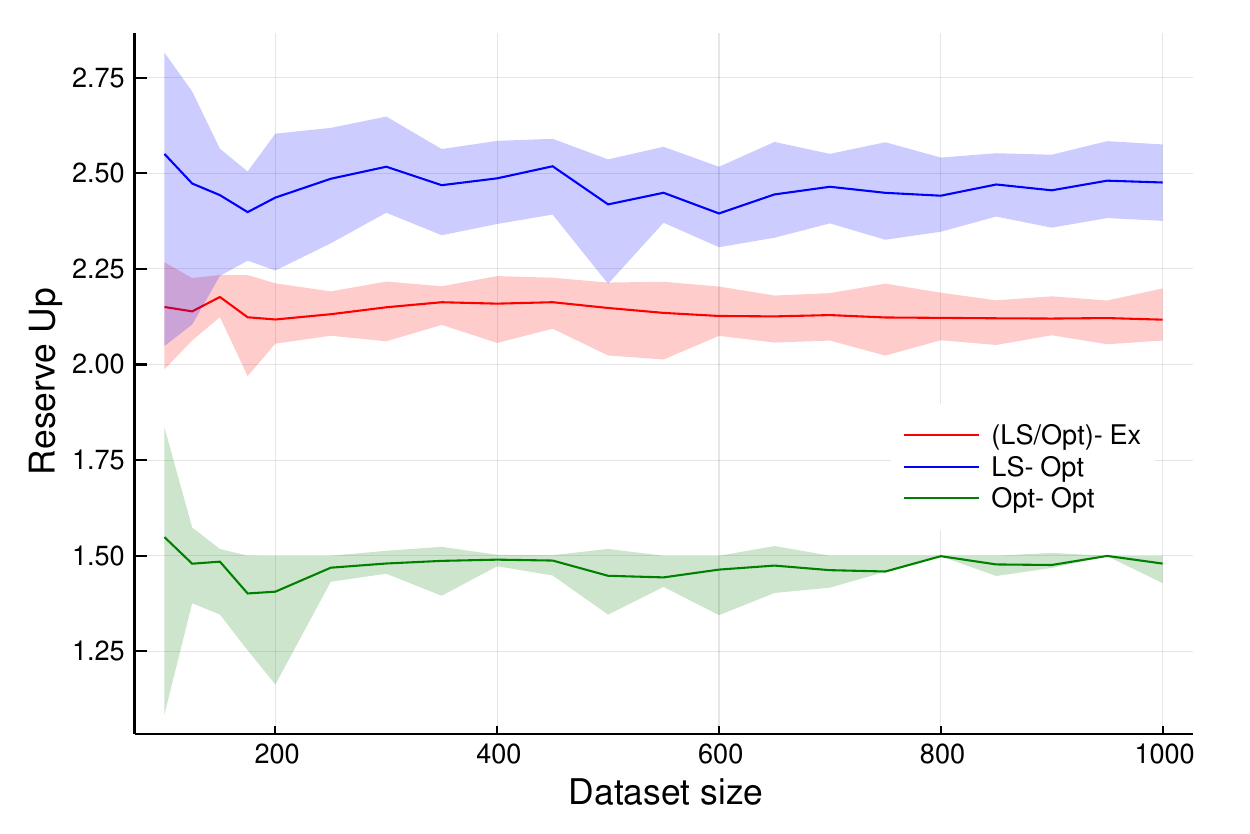}
}
\subfloat[Reserve down, $\theta_{R^{(dn)}}(0)$, from \eqref{mod-ar10-rdn}]{
\includegraphics[height=2.0in]{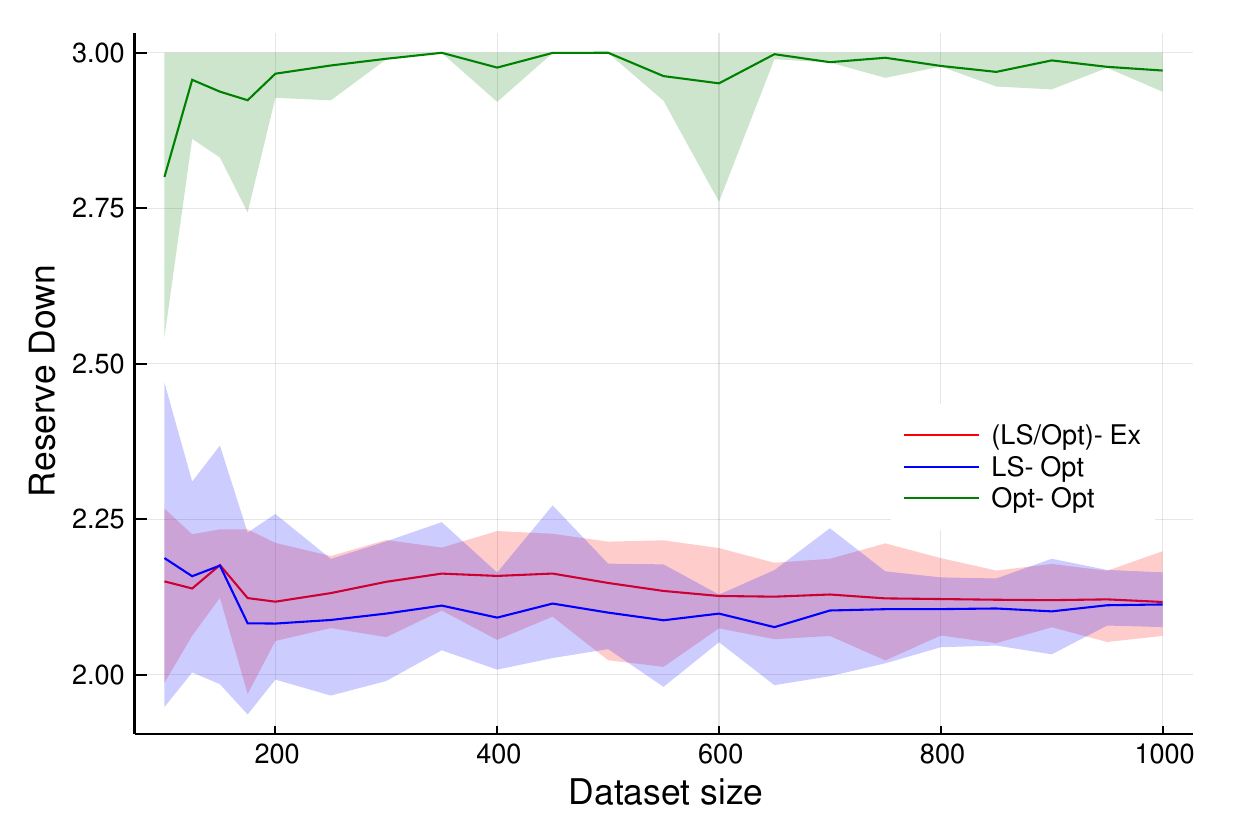}
}
\caption{Estimated parameters versus dataset size. Lines represent the average of the 100 estimation trials. Shaded areas represent the $10\%$ and $90\%$ quantiles. (a) and (b) Load coefficients, the models LS--Ex and LS--Opt coincide, thereby are presented as LS--(Ex/Opt). (c) and (d) Reserve coefficients The models Ls--Ex and Opt--Ex coincide, thereby are presented as (LS/Opt)--Ex.} \label{fig-demand coef-case1}
\end{figure}

In Figure \ref{fig-demand coef-case1} (c) and (d), we see that the Opt-Opt model greatly increases the downward reserve and decreases the upward reserve, both consistent with the change in the demand forecast parameters. Closed-loop estimation of only reserves led to increased up reserves that are the most expensive to violate, while downward reserves are mostly unaffected, this might be an artifact of the estimation model that uses the open-loop estimation as a starting point. The Opt-Opt model is limited to $3$ because that is the maximum reserve that can be allocated ($30\%$ of the generators' capacity).

To highlight the bias on load forecast, we present, in Figure \ref{fig-hist} (a), a histogram of deviations:
$error := realization - forecast$.
Negative values mean that the forecast value was above the realization. The LS estimation leads to an unbiased estimator, seen in the red histogram centered on zero. On the other hand, the forecast from the fully endogenous model is clearly biased, as it consistently forecasts higher values than the realizations. This fact is corroborated by the cumulative distribution functions displayed in Figure \ref{fig-hist} (b).

\begin{figure}[h!]
\subfloat[][\centering Two histograms are shown, the third color is their intersection (LS--Ex and LS--Opt coincide).]{
\includegraphics[height=2.0in]{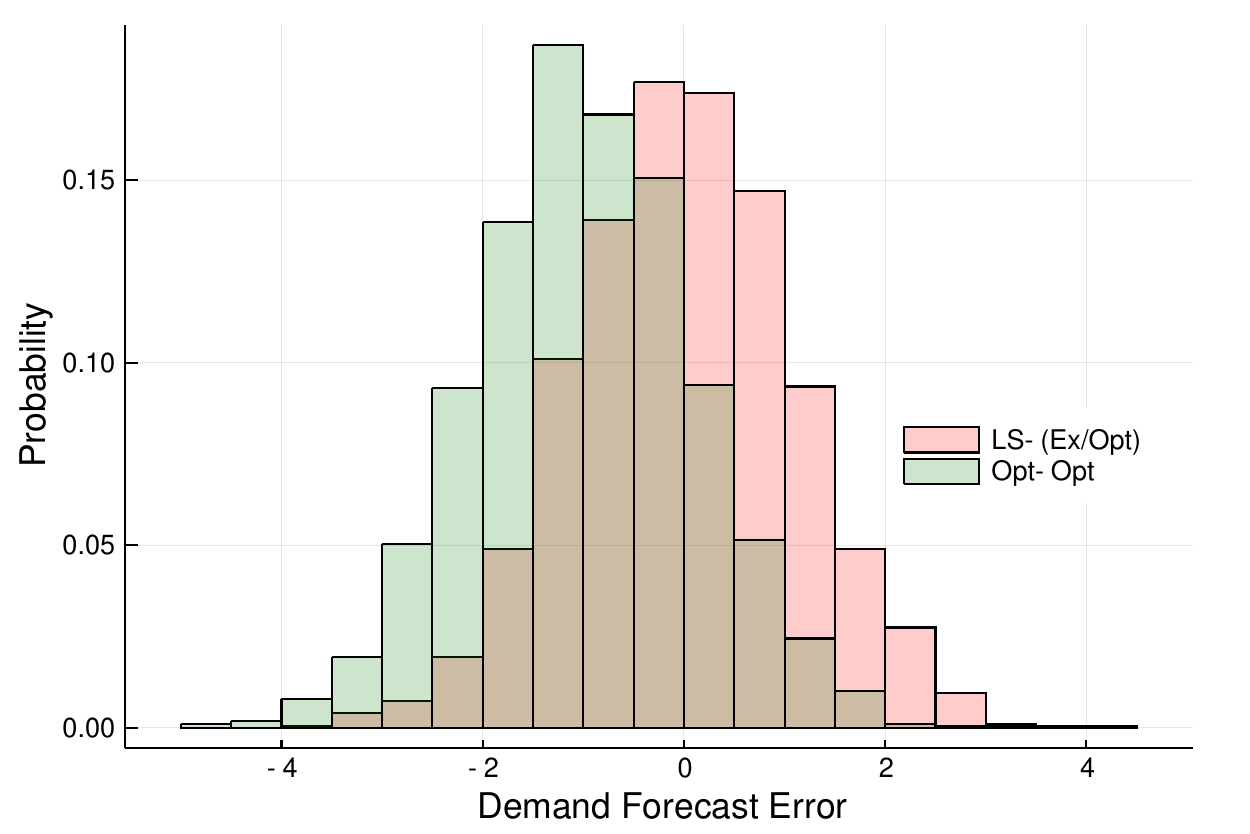}
}
\subfloat[][\centering Accumulated Probability -- out of the four models, three are shown here (LS--Ex and LS--Opt coincide).]{
\includegraphics[height=2.0in]{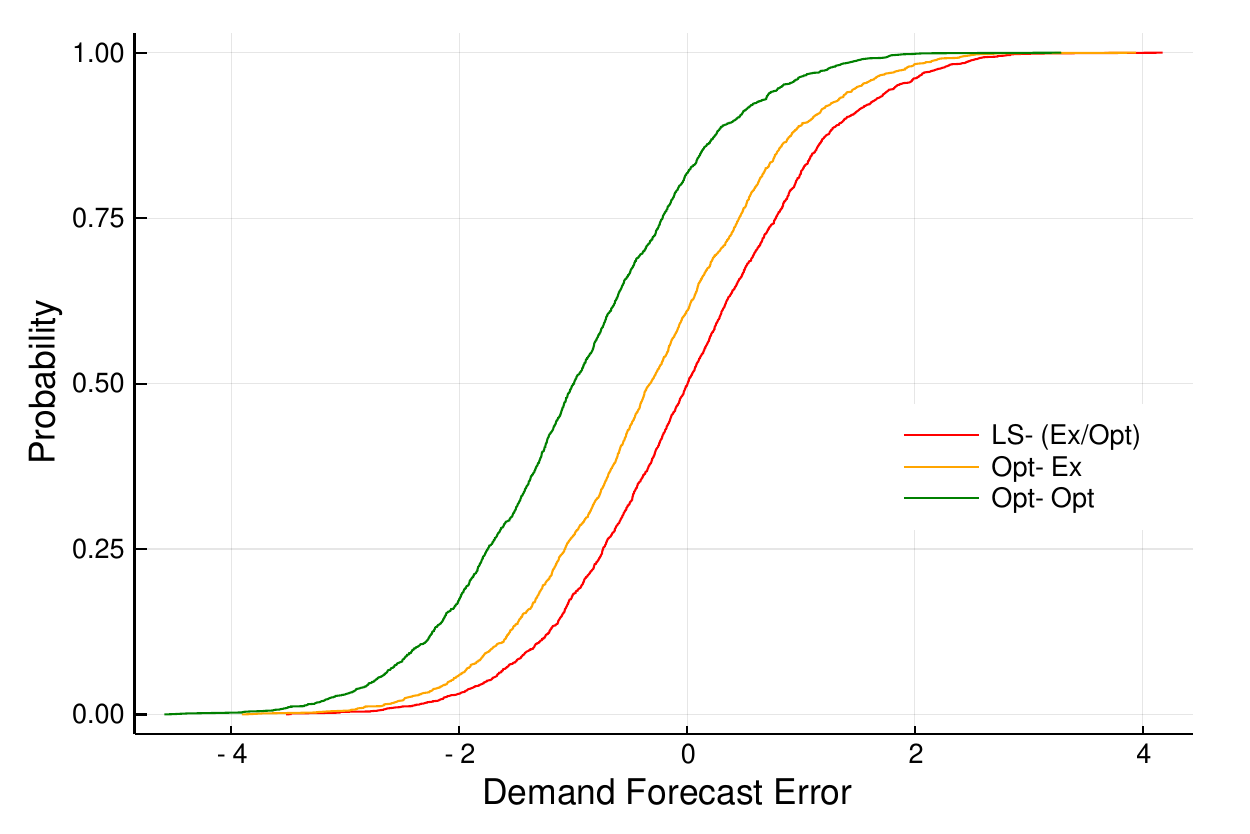}
}
\caption{Forecast error (observation -- forecast) in a histogram, comparing the fully optimized model with least squares estimation. Negative values mean the forecast was larger than the actual realization.}
\label{fig-hist} 
\end{figure}

\subsection{{\color{r2}Tests on large scale systems}}\label{secVeryLarge}

{\color{r2}
In this section, we emulate a realistic usage of the method, by 
testing it in realistic large-scale instances. The ones marked with ``\textit{rte}'' represent the French grid in 2013, using large-scale instances based on real systems \citep{josz2016ac}. The ones marked with ``\textit{pegase}'' are synthetic but based on the European grid \citep{fliscounakis2013contingency}.  All these instances are available in \cite{babaeinejadsarookolaee2019power}.
We limited the computing time to $30$ minutes so that the method can be used for dynamic reserve sizing and load forecast in hour-ahead dispatch.
Note that using the method in other time scales, such as day- or week-ahead dispatch, would also have great potential to benefit real systems around the globe.
The number of observations used for training was $600$, and the number of out-of-sample observations for testing was $10000$.
We used larger servers with $64$ processes in parallel and $1024$ Gb of RAM to perform the computations. These servers are available in Amazon Web Services, AWS, \citep{aws2022}.

\begingroup
\renewcommand{\arraystretch}{0.7} % Default value: 1

\begin{table}[h!]
\centering
{
\color{r2}
\begin{tabular}{l|rrrrr|rrrrr}
   & \multicolumn{5}{c}{Test} \vline & \multicolumn{5}{c}{Train} \\
\hline
Dataset & \multicolumn{2}{c}{Opt-Opt} & \multicolumn{2}{c}{LS-Opt} & LS-Ex & \multicolumn{2}{c}{Opt-Opt} & \multicolumn{2}{c}{LS-Opt} & LS-Ex \\
      & (\$) & (\%)  & (\$) & (\%) & (\$) & (\$) & (\%)  & (\$) & (\%)   & (\$)   \\
\hline
6468-rte     &     19870 &  8.56  &    20026 &  7.85  &  21731 &    20278 &  8.50  &   20462 &  7.67  &  22163  \\
6470-rte     &     25579 &  9.21  &    25849 &  8.25  &  28174 &    24004 & 10.38  &   24235 &  9.52  &  26785  \\
6495-rte     &     28726 &  9.77  &    28988 &  8.94  &  31835 &    28400 & 10.52  &   28666 &  9.68  &  31739  \\
6515-rte     &     34397 &  2.01  &    34557 &  1.56  &  35103 &    33757 &  2.01  &   33874 &  1.55  &  34354  \\
9241-pegase  &     76438 & 12.89  &    76662 & 12.62  &  87750 &    75315 & 10.31  &   75509 & 10.08  &  83973  \\
13659-pegase &     82608 &  7.97  &    82908 &  7.63  &  89760 &    83930 &  7.50  &   84196 &  7.21  &  90737

\end{tabular}
}
\caption{\color{r2}Results for very large-scale systems. LS-Ex is the reference, only costs in \$ are shown. For both Opt-Opt and LS-Opt costs are shown in \$, and their improvement compared to LS-Ex is shown in \%.}
\label{tabVeryLarge}
\end{table}

\endgroup

Results are shown in Table \ref{tabVeryLarge}, where the number in the dataset name represents the number of buses in the system. {\color{r4} The Opt-Opt model led to gains in the range of 2.01\% to 12.89\% (with an average equal to 8.42\%). These} results provide strong evidence that our method is capable of producing meaningful gains in practice.
}

\subsection{{\color{r3}Comparison with the least-squares and linear-bias benchmark approaches}}
\label{sec-realdata2}

{\color{r3}
Finally, in this last case study, we use the realistic system 6470-rte (from the previous Section) and real hourly load data from the PJM power system obtained from the U.S. Energy Information Agency \citep{eia2023}. The load profiles have a clear daily pattern. Hence, we expanded the previously used AR(1) to an AR(24) model for the load forecast. Each experiment comprises i) training the model in a sampled week from the historical hourly data, i.e., 168 observations were used in the training step, and ii) testing the estimated model out-of-sample in the following week, i.e., in the subsequent 168 observations from the historical data considered in the training procedure. For instance, if the week of a given training step ends at midnight June 3rd, the week for the test step starts at 1 AM June 4th. We also limited the training time to 30 minutes and used the same AWS servers with 64 cores and 1024 Gb of RAM to perform computations. We repeated the training-and-testing analysis for 10 different pairs of consecutive weeks to test the consistency of the gains.

Besides the typical least-squares benchmark (LS-Ex) and the proposed models (LS-Opt and Opt-Opt), we consider {\color{r4}the linear biasing approach from Section \ref{sec-naive} combined with the LS-Ex, labeled LS-Ex-Linear-Bias. To derive the linearly biased forecasts, we use the LS forecast as a reference and perform a grid search varying $\alpha$ from 1.0000 to 1.0500 in steps of 0.0025.}
Hence, $\alpha = 1.0000$ corresponds to the original LS-Ex model.

The results are summarized in Figure \ref{fig-realsys-realload}, which presents the out-of-sample percentage improvements of each model (LS-Ex-Linear-Bias, LS-Opt, and Opt-Opt) against the industry benchmark LS-Ex. We can see that Opt-Opt is consistently the best model, and LS-Opt is consistently the second best.

\begin{figure}[h!]
\begin{center}
\includegraphics[height=2.2in]{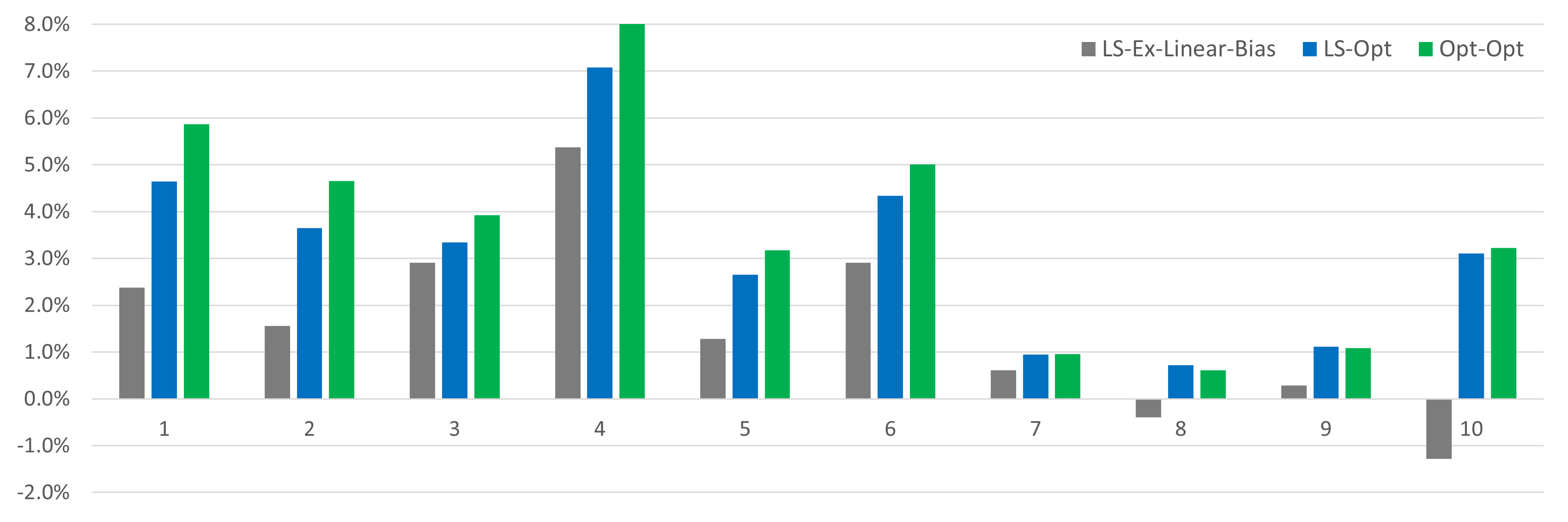}
\caption{Out-of-sample percentage improvement of proposed models compared to the LS-Ex benchmark. Improvement is shown for 10 experiments, each comprised of training the models in a week of data and testing the models in the following historical week.} \label{fig-realsys-realload}
\end{center}
\end{figure}

The best multipliers for the {\color{r4}LS-Ex-Linear-Bias} ranged between 1.0000 and 1.0100 (the most frequent was 1.0050) for this power system, which led to improvements in 8 out of 10 experiments against the LS-Ex benchmark. The losses observed in two cases are expected in similar procedures where a model is selected based on training data and then deployed on unseen data. Some consecutive weeks might have a significant change in the load behavior. Despite that, this simple {\color{r4}linear biasing strategy} showcases a relevant benefit, which, in practice, can be further improved by means of different reestimation policies that, for instance, increase the re-estimation frequency and the relevance of the most recent observations. 

We also notice a reduced improvement of both Opt-Opt and LS-Opt in experiments 7, 8, and 9. We attribute this to a higher average load in these weeks. The higher load makes the system dispatch expensive resources, which leads to higher overall costs and decreases the percentage gain. Moreover, in such situations, the system is so constrained that changing the forecast and the reserve margins cannot significantly modify the pre- and post-dispatch setting points in a substantial way.

This final case study brings strong evidence corroborating our thesis that Opt-Opt and LS-Opt models can be used in practice and lead to meaningful gains in real power systems with thousands of buses and nontrivial load profiles.}

\section{Conclusions}

We have proposed a mathematical framework for application-driven learning, a technique whereby forecast errors are measured in terms of the underlying objective function of the problem. While similar ideas have been developed in recent literature under multiple names, such as end-to-end learning, or integrated estimation and optimization, for example, our approach allows for the study of the case where the underlying problem is a two-stage model for which one is interested in solving it using a point forecast of the uncertainty. We do so by building upon the ideas of bilevel optimization. 

{\color{r4}Motivated by the electricity market practice, wherein relevant evidence points toward system operators deviating their demand forecasts from the minimum error estimates (e.g., least-squares) to obtain better operational results, we apply our new learning framework to find the best load and reserve requirement forecasting model for the application of energy and reserve scheduling.}

{\color{r4}Numerical experiments 
demonstrate the effectiveness of the proposed closed-loop approach in minimizing the system cost in the long run. For the large systems, we found gains in the range of 2.01\% to 12.89\% (with an average equal to 8.42\%), thereby, providing strong evidence that our method is capable of producing relevant gains in practice. Moreover, we have empirically shown that the method scales well and is able to handle large-scale realistic power systems---in the largest example, 13,659 buses with real data---with a limited training time of $30$ minutes resorting to cloud computing. Consequently, the method might be well fit for realistic hour-ahead planning, in addition to day- and week-ahead.}

The proposed framework is fairly general, but we focused on right-hand side uncertainty and linear models as this was the relevant setting for the reserve requirement forecasting problem. Possible extensions of this work could aim at generalizing this setting. Indeed, the exact method would work for objective-function uncertainty and for non-linear models with strong duality, like many conic programs. The heuristic method is even more flexible, it could accommodate uncertainty anywhere in the models, and it could be used in many non-linear (even integer) problems. {\color{r2}Also, the heuristic method can benefit from online optimization approaches and warm starts based on the solution of a previous execution of the algorithm with similar data as would occur in a real-time application. This could bring relevant gains to the optimization of sequential hour-ahead operations.} 
Other aspects of estimation procedures to be explored include regularization with shrinkage operators \citep{tibshirani2011regression},  or the addition of LS estimates in the objective \citep{kao2009directed} with a penalizing coefficient to provide a balance between classical and application-driven forecasts. Moreover, it would be interesting to conduct experiments to verify the empirical performance of various statistical models ($\psi$), such as vector autoregressive models {\color{r4}and neural networks}. Proving convergence conditions {\color{r4}and rates} for many of the aforementioned changes in the model would provide significant contributions to the literature.

\section*{Acknowledgments}

The work of Joaquim Dias Garcia was partially supported by the Coordenação de Aperfeiçoamento de Pessoal de Nível Superior - Brasil (CAPES) - Finance Code 001. 
The work of Alexandre Street was partially supported by Fundação de Amparo à Pesquisa do Estado do Rio de Janeiro (FAPERJ) and Conselho Nacional de Desenvolvimento Científico e Tecnológico (CNPq).
Tito Homem-de-Mello acknowledges the support of grant FONDECYT 1221770  from ANID, Chile. 
The authors also thank PSR for making PSRCloud available for the experiments in Sections \ref{secVeryLarge} and \ref{sec-realdata2}.

\begin{APPENDICES}

\section{Proofs of Theorems}\label{sec:proof}

\begin{proof}{Proof of Theorem \ref{thm}:}
First, notice that $G_i(z,Y)$, $i \in \{a, p\}$, is continuous with respect to its arguments as it is a sum of a linear function and the optimal value of a parametric program \citep{gal2010postoptimal}.
Recall that $\zeta$ is a continuous vector-valued function because of Assumption \ref{as1}.
Hence, $G_a(\zeta(\Psi(\theta,X)),Y)$ is a real-valued continuous function. Next, we show that $G_a(\zeta(\Psi(\theta,X)),Y)$ is integrable.
Indeed, since $Z$ is bounded (Assumption \ref{as_z}), it follows that $\zeta(y)$ is bounded for all $x$ by a constant, say $K_1$, so that $\|\zeta(y)\| \leq K_1$.
By duality, $Q_i(z,y) = \max \limits_{\pi} \{(b_i - H_i z + F_i y)^{\top} \pi \ | \ W_i^{\top} \pi = q_i, \pi \geq 0 \}$, but by Assumption \ref{as_dualQ} the dual variables of $Q_a(z,y)$ are bounded by a constant, say $K_2$, so $\|\pi\| \leq K_2$. Thus, by a sequence of applications of Cauchy-Schwarz and triangle inequalities, we have that
\begin{align*}
  \big|G_a(\zeta(\Psi(\theta,X)),Y)\big| &\leq \big|c_a^{\top}\zeta(\Psi(\theta,X)) + Q_a(\zeta(\Psi(\theta,X)),Y)\big| \\
  & \leq \big\|c_a\big\| \big\|\zeta(\Psi(\theta,X))\big\| + \big\|b_a - H_a \zeta(\Psi(\theta,X)) + F_a Y\big\| \max \limits_{ W_i^{\top} \pi = q_a, \pi \geq 0}\big\|\pi\big\| \\
  & \leq K_1 \big\|c_a\big\| + K_2 \left(\big\|b_a\big\| + \big\|H_a \zeta(\Psi(\theta,X))\big\| + \big\| F_a Y\big\|\right)\\ 
  & \leq K_1 \big\|c_a\big\| + K_2 \left(\big\|b_a\big\| + \big\|H_a\big\|  K_1 + \big\| F_a \big\| \big\|Y\big\|\right). 
\end{align*}
Hence, since $Y$ is integrable (condition (iv) of the Theorem), we have that  $G_a(\zeta(\Psi(\theta,X)),Y)$ is integrable.

It follows that the conditions of Theorem
7.53 in \citet{Shapiro:09} are satisfied and we conclude that: (i)
the function $\varphi(\theta):=\mathbb{E}[G_a(\zeta(\Psi(\theta,X)),Y)]$
is finite valued and continuous in $\theta$,  (ii) by the Strong Law of
Large Numbers, for any $\theta\in\Theta$ we have
\begin{align}
\lim_{T\to\infty}\frac{1}{T}\sum_{t=1}^{T} & G_a\big(\zeta(\Psi(\theta, X_t)), Y_t\big)=\Ea{G_a\big(\zeta(\Psi(\theta,X)),Y\big)}\ \ \text{w.p.1},\label{eq:p:convergence}
\end{align} and (iii) the convergence
in \eqref{eq:p:convergence} is \emph{uniform} in $\theta$. Thus, by Theorem 5.3
in \citet{Shapiro:09}, since the set $\Theta$ is compact, we have that the minimizers (over $\Theta$) of the expression {\color{r2}inside the limit} on the left-hand side of \eqref{eq:p:convergence}---i.e., $\theta_T$---converge to the minimizers of the expression on the right-hand side in the sense of  \eqref{eq-conv-set}--\eqref{eq-conv-obj}. \textit{Q.E.D.}
\end{proof}

\vspace{1cm}

For reference, we state now a lemma that provides a more general result than Theorem 7.53
in \citet{Shapiro:09}. {\color{r3}We remark that the result in the lemma is not new; in fact, \cite{korf2001random} show \textit{epi-convergence} of long-run averages of functions of ergodic processes, which implies uniform convergence. \cite{surbirge:21} also show a similar result for Harris recurrent stationary Markov chains.  These results notwithstanding, we present the lemma for the sake of completeness since its proof is a straightforward extension of the proof of Theorem 7.53 in \citet{Shapiro:09}. }

\begin{lem}
\label{lemma:SDR_ergodic}
Theorem 7.53 in \citet{Shapiro:09} is still valid if the i.i.d.\ assumption is replaced with a weaker assumption that the samples form a 
stationary ergodic process.
\end{lem}
\begin{proof}{Proof:}
  Any measurable function of a stationary ergodic process is also a stationary ergodic process \citep{billingsley1986}. Moreover, if a process  $\{W_t\}_{t=1}^{\infty}$ is stationary and ergodic, then the classical ergodic  theorem (see, e.g., \cite{billingsley1986}) ensures that 
\[
\lim_{T\to\infty}\frac{1}{T}\sum_{t=1}^{T}  W_t \ =\ \Ea{W_1}\ \ \text{w.p.1}.
\] 
A closer look at the proof of Theorem 7.53 in \citet{Shapiro:09} shows that the i.i.d.\ assumption is used only to invoke the Strong Law of Large Numbers, which, as shown above, can be replaced by the ergodic theorem in the more general case.   \textit{Q.E.D.}
\end{proof}

\vspace{1cm}

\begin{proof}{Proof of Theorem \ref{thm2}:}
Fix $\theta \in \Theta$. Consider the function 
$\Phi$ defined as 
\[
\Phi(Y_1,\ldots,Y_t)\ := \ G_a\big(\zeta(\Psi(\theta, X_t)), Y_t\big)
\]
and process $\{W_t\}_{t=1}^{\infty}$ defined as $W_t:=\Phi(Y_1,\ldots,Y_t)$. Under the assumption that $\{ Y_t\}_{t=1}^{\infty}$ is a stationary ergodic time series, it follows that $\{W_t\}_{t=1}^{\infty}$ is stationary and ergodic, since $\Phi$ is measurable. Lemma~\ref{lemma:SDR_ergodic} 
shows that  \eqref{eq:p:convergence} holds in this case and, therefore, the proof follows the same steps as those of the proof of Theorem~\ref{thm}.  \textit{Q.E.D.}
\end{proof}

\vspace{1cm}

\begin{proof}{Proof of Corollary \ref{cor}:}
Because of Theorem~\ref{thm} or Theorem~\ref{thm2}, we know that the distance between $\theta_{T}$ and the set $S^*$ of minimizers of the underlying problem (given by \eqref{eq-conv-obj}) converges to zero as $T\to\infty$. Thus, any convergent subsequence of the process $\{\theta_T\}_{T=1}^{\infty}$ converges to some minimizer $\theta^*\in S^*$. Moreover, since $\{\theta_T\}_{T=1}^{\infty}\subseteq \Theta$ and $\Theta$ is assumed to be compact, it follows that there exists at least one convergent subsequence. Meanwhile, any other possible choice of $\theta$, obtained by any other estimation method, is merely a feasible solution for the underlying problem and, hence, cannot be strictly better than $\theta^*$. Q.E.D.
\end{proof}

\section{Complete bilevel formulation of Application-Driven Load Forecasting and Reserve Sizing}\label{app-full}
\begin{align}
& \min_{\substack{\theta_D, \theta_{R^{up}}, \theta_{R^{dn}}, \\ \hat{D}_t, \hat{R}_t^{(up)}, \hat{R}_t^{(dn)}, g_t, \delta_t^{LS}, \delta_t^{SP} \\ g_t^* , r_t^{(up)*}, r_t^{(dn)*}}} 
        \frac{1}{T}\sum_{t \in \mathbb{T}}
            \big[c^{\top} g_t^* + p^{(up)\top} \hat{r_t}^{(up)*} + p^{(dn)\top} \hat{r_t}^{(dn)*} + 
            \lambda^{LS}\delta_t^{LS} + \lambda^{SP}\delta_t^{SP} \big] \hspace{1.0cm} \label{Fmod-pow-op-obj}\\
& \hspace{15px} s.t. \hspace{18px} \forall t \in \mathbb{T}: \notag \\
& \hspace{52px} \hat{D}_t = \Psi_D (\theta_D, x_t) \label{Fmod-pow-fore-d}\\
& \hspace{52px}   \hat{R}_t^{(up)} = \Psi_{R^{(up)}} (\theta_{R^{(up)}}, x_t) \label{Fmod-pow-fore-up}\\
& \hspace{52px}   \hat{R}_t^{(dn)} = \Psi_{R^{(dn)}} (\theta_{R^{(dn)}}, x_t) \label{Fmod-pow-fore-dn}\\
& \hspace{52px} e^{\top} (M g_t - \delta_t^{SP}) = e^\top (D_t - \delta_t^{LS}) \label{Fmod-pow-op-bal}\\
& \hspace{52px} -F \leq B(M g_t + \delta_t^{LS} - D_t - \delta_t^{SP}) \leq F  \label{Fmod-pow-op-net}\\
& \hspace{52px} g_t^* - r_t^{(dn)*}\leq g_t \leq g_t^* + r_t^{(up)*}  \label{Fmod-pow-op-adj}\\
& \hspace{52px} \delta_t^{LS}, \delta_t^{SP}, \hat{R}_t^{(up)}, \hat{R}_t^{(dn)}, g_t \geq 0 \label{Fmod-pow-op-zeeo}\\
% TODO fix this long line
& \hspace{42px} \Big( g_t^* , r_t^{(up)*}, r_t^{(dn)*}\Big) \in  \argmin_{\substack{\hat{g_t}, \hat{\delta}_t^{LS}, \hat{\delta}_t^{SP},\\\hat{r_t}^{(up)}, \hat{r_t}^{(dn)}}}  \big[\tilde{c}^{\top} \hat{g_t} + \tilde{p}^{(up)\top} \hat{r_t}^{(up)} + \tilde{p}^{(dn)\top} \hat{r_t}^{(dn)} +  
            \tilde{\lambda}^{LS}\hat{\delta}_t^{LS} + \tilde{\lambda}^{SP}\hat{\delta}_t^{SP} \big]  \label{Fmod-pow-pl-obj}\\
& \hspace{155px} s.t. \hspace{30px} e^\top (M \hat{g_t} - \hat{\delta}_t^{SP}) = e^\top (\hat{D}_t - \hat{\delta}_t^{LS}) \label{Fmod-pow-pl-bal}\\
& \hspace{204px} -\tilde{F} \leq \tilde{B}(M\hat{g_t} + \hat{\delta}_t^{LS} - \hat{D}_t - \hat{\delta}_t^{SP}) \leq \tilde{F} \label{Fmod-pow-pl-net}\\
& \hspace{204px}  N \hat{r_t}^{(up)} = \hat{R}_t^{(up)} \label{Fmod-pow-pl-balup}\\
& \hspace{204px} N \hat{r_t}^{(dn)} = \hat{R}_t^{(dn)} \label{Fmod-pow-pl-baldn}\\
& \hspace{204px} \hat{g_t} + \hat{r_t}^{(up)} \leq \tilde{K} \label{Fmod-pow-pl-ub}\\
& \hspace{204px} \hat{g_t} - \hat{r_t}^{(dn)} \geq 0 \label{Fmod-pow-pl-lb}\\
& \hspace{204px} \hat{r_t}^{(up)} \leq {\bar{r}}^{(up)} \label{Fmod-pow-pl-upub}\\
& \hspace{204px} \hat{r_t}^{(dn)} \leq {\bar{r}}^{(dn)} \label{Fmod-pow-pl-dnub}\\
& \hspace{204px} \hat{g_t}, \hat{r_t}^{(up)}, \hat{r_t}^{(dn)}, \hat{\delta}_t^{LS}, \hat{\delta}_t^{SP} \geq 0. \label{Fmod-pow-pl-zero}
\end{align}

\section{Scalability of the proposed method}\label{sec:scale}

We analyze the algorithm's scalability and its consequent applicability to a broad range of large power system networks. We created instances with the number of buses ranging from $600$ to $6000$. These instances were created by connecting multiple copies of the IEEE $300$ bus case from PG-LIB-OPF \citep{babaeinejadsarookolaee2019power} (with flow limits reduced by $25\%$ and loads multiplied by $0.9$, to stress the system). The optimization was performed on an Intel Xeon E5-2680 with 12 cores at 2.50GHz, 128Gb RAM. We generated a single training dataset with $1000$ observations for each instance and optimized each problem considering four- and twelve-hour time limits. Then, we evaluated the solutions obtained with each method and computational time limit with a common dataset of $10000$ out-of-sample observations.

The results, in terms of Test and Train cost, are depicted in Table \ref{tabLarge}. It is clear that the proposed methods (LS-Opt and Opt-Opt) consistently outperformed the least-squares benchmark (LS-Ex) despite reaching the time limit for training. Moreover, the co-optimization scheme (Opt-Opt) was better than the other two methods (LS-Ex and LS-Opt) in both in-sample and out-of-sample evaluation. Also, the method performed well, generalizing from training to testing. For the larger systems, we see smaller relative improvements in the cost functions. This is due to (i) the increased dimension of linear programs solved in each iteration, and (ii) the number of parameters, which increases with the system size. Reasons (i) and (ii)  imply that, for a given computational budget, fewer iterations are run in the training stage, possibly yielding sub-optimal solutions. Thus, it is conceivable that much better results could be obtained with more time (or processing capacity), as the $12$-hour runs lead to improvements that are more than three times those obtained with $4$-hour runs. Although we were limited to $12$ cores in this experiment, one could use one core (or more) per observation in the training dataset, leading to massive speed-ups. 
It is worth mentioning that, in all cases, the LS estimation found coefficients very close to the true ones. This fact, together with the gains shown in  Table \ref{tabLarge} for many instances of different sizes, indicates that the forecast bias introduced by our methodology is consistent in promoting an improved operation. Notably, relevant gains were found even in cases of very large-scale power systems. For instance, for the 3000-bus system, in this simulation run, the gain was $13\%$ with the $12$-hour time limit. Therefore, the results provide strong evidence that our method is capable of producing meaningful gains in practice.
\begingroup
\renewcommand{\arraystretch}{0.7} % Default value: 1
\begin{table}[h!]
\centering
\begin{tabular}{rr|rrrrr|rrrrr}
  &  & \multicolumn{5}{c}{Test}       \vline    & \multicolumn{5}{c}{Train}                     \\
\hline
Buses & Time    & \multicolumn{2}{c}{Opt-Opt}   & \multicolumn{2}{c}{LS-Opt}    & LS-Ex    & \multicolumn{2}{c}{Opt-Opt}  & \multicolumn{2}{c}{LS-Opt}   & LS-Ex    \\
     & (h) & (\$) & (\%)  & (\$) & (\%) & (\$) & (\$) & (\%)  & (\$) & (\%)   & (\$)   \\
\hline
600  &  4  &     15256 &  21.51  &    16741 & 13.87  &  19437 &    14898 &  22.73  &   16552 & 14.16  &  19281 \\
600  & 12  &     15226 &  21.66  &    16739 & 13.88  &  19437 &    14873 &  22.86  &   16563 & 14.10  &  19281 \\
1200 &  4  &     35261 &  26.38  &    40103 & 16.28  &  47899 &    34985 &  26.87  &   39994 & 16.41  &  47843 \\
1200 & 12  &     35158 &  26.60  &    38843 & 18.91  &  47899 &    34842 &  27.17  &   38787 & 18.93  &  47843 \\
1800 &  4  &     60043 &  18.61  &    69355 &  5.98  &  73770 &    61033 &  17.97  &   69915 &  6.03  &  74405 \\
1800 & 12  &     53378 &  27.64  &    62786 & 14.89  &  73770 &    54294 &  27.03  &   63224 & 15.03  &  74405 \\
2400 &  4  &     81018 &  20.23  &    94675 &  6.78  & 101561 &    78686 &  19.57  &   91278 &  6.70  &  97832 \\
2400 & 12  &     80330 &  20.91  &    93739 &  7.70  & 101561 &    77969 &  20.30  &   90413 &  7.58  &  97832 \\
3000 &  4  &    120389 &   5.18  &   125465 &  1.19  & 126972 &   119214 &   5.27  &  124328 &  1.20  & 125841 \\
3000 & 12  &    110302 &  13.13  &   122560 &  3.47  & 126972 &   109231 &  13.20  &  121382 &  3.54  & 125841 \\
3600 &  4  &    149141 &   3.51  &   153484 &  0.69  & 154558 &   147110 &   3.50  &  151404 &  0.68  & 152439 \\
3600 & 12  &    136479 &  11.70  &   150303 &  2.75  & 154558 &   134664 &  11.66  &  148332 &  2.69  & 152439 \\
4200 &  4  &    177451 &   1.72  &   179785 &  0.43  & 180555 &   174395 &   1.70  &  176651 &  0.43  & 177406 \\
4200 & 12  &    165963 &   8.08  &   177444 &  1.72  & 180555 &   162820 &   8.22  &  174309 &  1.75  & 177406 \\
4800 &  4  &    206358 &   1.22  &   208300 &  0.29  & 208910 &   209816 &   1.28  &  211927 &  0.29  & 212546 \\
4800 & 12  &    197707 &   5.36  &   206539 &  1.13  & 208910 &   201294 &   5.29  &  210197 &  1.11  & 212546 \\
5400 &  4  &    232071 &   1.05  &   233992 &  0.23  & 234524 &   228660 &   1.23  &  231018 &  0.21  & 231514 \\
5400 & 12  &    225203 &   3.97  &   232658 &  0.80  & 234524 &   222018 &   4.10  &  229769 &  0.75  & 231514 \\
6000 &  4  &    260427 &   0.96  &   262493 &  0.18  & 262963 &   262049 &   0.99  &  264198 &  0.18  & 264668 \\
6000 & 12  &    255384 &   2.88  &   261500 &  0.56  & 262963 &   257062 &   2.87  &  263118 &  0.59  & 264668
\end{tabular}
\caption{Results for very large-scale systems. LS-Ex is the reference, only costs in \$ are shown. For both Opt-Opt and LS-Opt costs are shown in \$ and their improvement compared to LS-Ex is shown in \%.}
\label{tabLarge}
% \vspace{-4mm}
\end{table}
\endgroup

\end{APPENDICES}

%%%%%%%%%%%%%%%%%%%%%%%%%%%%%%%%%%%%%%%%%%%%%%%%%%%%%%%%%%%%%%%
% References here (outcomment the appropriate case)

% CASE 1: BiBTeX used to constantly update the references
%   (while the paper is being written).
\bibliographystyle{informs2014} % outcomment this and next line in Case 1
\bibliography{refs} % if more than one, comma separated

% CASE 2: BiBTeX used to generate mypaper.bbl (to be further fine tuned)
%\input{mypaper.bbl} % outcomment this line in Case 2

%If you don't use BiBTex, you can manually itemize references as shown below.

%\begin{thebibliography}{}
%\bibitem[{American Butter Institute(2005)}]{abi}
%American Butter Institute (2005) Dairy market report. Retrieved June
%14, 2005, www.butterinstitute.org.
%\end{thebibliography}
%%%%%%%%%%%%%%%%%%%%%%%%%%%%%%%%%%%%%%%%%%%%%%%%%%%%%%%%%%%%%%

%%%%%%%%%%%%%%%%%
\end{document}